\documentclass{amsart}

\usepackage[latin1]{inputenc}

\usepackage{amsmath,amssymb,amscd,latexsym}

\usepackage[dvips]{graphicx}
\usepackage[english]{babel}

\def\Den{S}
\newcommand{\Snu}{\mathbb{S}_{\nu}}
\def\Coex{\text{Coex}}
\def\SepCoex{\text{Sep-Coex}}
\def\Act{\text{Act}}

\def\blue{\text{blue}}
\def\yellow{\text{yellow}}
\def\green{\text{green}}
\def\bluestar{\text{blue}^*}
\def\yellowstar{\text{yellow}^*}
\def\greenstar{\text{green}^*}
\def\sy{s_{\text{yellow}}}
\def\sb{s_{\text{blue}}}

\newcommand{\bor}[1][]{\mathcal{B}(#1)}

\newcommand{\communique}{\leftrightarrow}
\newcommand{\N}{\mathbb{Z}_{+}}

\newcommand{\Z}{\mathbb{Z}}

\newcommand{\Zd}{\mathbb{Z}^d}

\newcommand{\Pcond}{\overline{\mathbb{P}}}
\newcommand{\R}{\mathbb{R}}
\newcommand{\Rd}{\mathbb{R}^d}

\renewcommand{\P}{\mathbb{P}}
\newcommand{\E}{\mathbb{E}\ }
\newcommand{\F}{\E_{\Pcond}\ }

\newcommand{\Ed}{\mathbb{E}^d}

\newcommand{\Ber}{\text{Ber}}

\newcommand{\as}{\text{ a.s.}}

\renewcommand{\epsilon}{\varepsilon}
\renewcommand{\limsup}{\overline{\lim}}

\newcommand{\iid}{\text{i.i.d. }}
\newcommand{\ie}{\emph{i.e. }}
\newcommand{\resp}{\emph{resp. }}

\newcommand{\miniop}[3]{%
\renewcommand{\arraystretch}{0.6}
\begin{array}{c}
{\scriptstyle #1}\\
#2\\
{\scriptstyle #3}
\end{array}
\renewcommand{\arraystretch}{1}}

\newtheorem{theorem}{Theorem}[section]
\newtheorem{lemme}[theorem]{Lemma}
\newtheorem{defi}[theorem]{Definition}
\newtheorem{coro}[theorem]{Corollary}

\title[Coexistence in two-type models]{Coexistence in two-type first-passage percolation models}

\begin{document}

\author{Olivier Garet}
\author{R{\'e}gine Marchand}
\address{Laboratoire de Math{\'e}matiques, Applications et Physique
Math{\'e}matique d'Orl{\'e}ans UMR 6628\\ Universit{\'e} d'Orl{\'e}ans\\ B.P.
6759\\
 45067 Orl{\'e}ans Cedex 2 France}
\email{Olivier.Garet@labomath.univ-orleans.fr}

\address{Institut Elie Cartan Nancy (math{\'e}matiques)\\
Universit{\'e} Henri Poincar{\'e} Nancy 1\\
Campus Scientifique, BP 239 \\
54506 Vandoeuvre-l{\`e}s-Nancy  Cedex France}
\email{Regine.Marchand@iecn.u-nancy.fr}


\begin{abstract}
We study the problem of  coexistence in a two-type competition model governed
by first-passage percolation on $\Zd$ or on the infinite cluster in Bernoulli percolation.
Actually, we prove for a large class of ergodic stationary passage times that for distinct points $x,y\in\Zd$, there is a strictly positive probability that
$\{z\in\Zd;d(y,z)<d(x,z)\}$ and $\{z\in\Zd;d(y,z)>d(x,z)\}$ are both infinite sets.
We also show that there is a strictly positive probability that the graph of time-minimizing path from the origin in first-passage percolation has at least two topological ends.   
This generalizes results obtained by H{\"a}ggstr{\"o}m and Pemantle for
independent exponential times on the square lattice.
\end{abstract}

\maketitle

\setcounter{tocdepth}{1}     
\setcounter{secnumdepth}{1}  


\section{Introduction}


The two-type Richardson's model was introduced by H{\"a}ggstr{\"o}m and Pemantle in~\cite{Haggstrom-Pemantle-1} as a simple competition model between two infections: on the cubic grid $\Z^d$, two
distinct infections, type 1 and type 2, starting respectively from two
distinct sources $s_1,s_2 \in \Z^d$, compete to invade the sites of the
grid $\Z^d$. Each one progresses like a first-passage percolation process
on $\Z^d$, governed by the same  family $(t(e))_{e \in\Ed}$ of \iid 
exponential random variables, indexed by the set $\Ed$ of edges of $\Z^d$, but the two infections interfere in the
following way: once a site is infected by the type $i$ infection, it 
remains of type $i$ forever and can not transmit the other infection. This leads to two very different possible evolutions of the process:
\begin{itemize}
\item
either one infection surrounds the other one, stops it and then goes on infecting the remaining healthy sites as if it was alone, 
\item
or the two infections grow mutually unboundedly, which is called  \textit{coexistence}.
\end{itemize}
The probability that, given two distinct sources, coexistence occurs is of course not full, and the relevant question is to determine whether coexistence occurs with positive probability or not. Although this competition problem is interesting in its own right, it is also a powerful
tool to study the existence of two semi-infinite geodesics (or topological ends) of the embedded spanning tree in the related first-passage percolation model.
Thus, H{\"a}ggstr{\"o}m and Pemantle  proved that coexistence for any two initial sources in the two-type Richardson's model on $\Z^2$ occurs with positive probability, and consequently that in first-passage percolation on $\Z^2$ with \iid exponential passage times on the edges, the probability that there exist at least two topological ends  in the embedded spanning tree is positive.

Their results strongly rely  on an interacting particle representation of the problem which is typical of the exponential passage times. The aim of this paper is to extend these results to more general passage times, where this representation is not available anymore or at least much less natural.
We consider here stationary ergodic first-passage percolation on $\Zd$, $d\ge 2$  (and also on an infinite cluster of Bernoulli percolation), and  prove that, for any two distinct sources, the probability that coexistence occurs is strictly positive. As a consequence, we obtain that in the related first-passage percolation on $\Z^d$, the probability that there exist at least two topological ends  in the embedded spanning tree is positive.

The structure of the proof is the following. First, the key step is to prove that there exist two sources such that coexistence occurs, and this is the aim of Section~3. Heuristically, the shape theorem of first-passage percolation, combined with the fact that the two infections have the same speed, gives the intuition that the largest the distance between the two sources is, the hardest it is for one infection to surround the other one. More precisely, Theorem~\ref{yarfyarf} says that
if $d(x,y)$ denotes the 
travel time between the sites $x$ and $y$, then there exists a  site $x$ such that the event 
\begin{itemize}
\item
\textit{and } the set of sites $z$ such that $d(0,z)<d(x,z)$ is infinite,
\item
\textit{and } the set of sites $z$ such that $d(0,z)>d(x,z)$ is infinite
\end{itemize}
has positive probability.
The proof of this result relies on the existence of a directional  asymptotic speed in the
related first-passage percolation model. 

The next step is to transfer the coexistence result for these sources to \textit{any} two  initial sources; this is done by a modification argument of the configuration around the sources using a finite energy property for the passage times. Roughly speaking, this result expresses the fact that non-coexistence is due to a \textit{local} advantage obtained by one infection at the first moments of the competition. The two topological ends result is shown by a similar modification argument. These results are proved separately in Section~4 for diffuse passage times and in Section~5 for  integer passage times.

The last section is finally devoted to the study of a probabilistic cellular automata describing a discrete competition model between two infection types, related to the chemical distance in super-critical Bernoulli percolation on $\Z^d$.\\

We start now with a reminder of the result of existence of directional asymptotic speeds in classical first-passage percolation, and an extension of this result to  first-passage percolation on an infinite Bernoulli cluster. 

 
\section{Reminder on the directional asymptotic speed results}


In classical first-passage percolation, one has the well-known directional asymptotic speed result:  if $(t(e))_{e
  \in \Ed}$ are \iid non-negative integrable random variables, then for
every $x\in \Z^d$, there exists $\mu(x)\geq 0$ such that a.s :
$$\lim_{n \rightarrow \infty}\frac{t(0,nx)}{n}=\lim_{n \rightarrow
  \infty}\frac{\E t(0,nx)}{n}=\mu(x).$$

This result has been extended in full details in a previous work of the authors \cite{garet-marchand} to first-passage percolation on an infinite Bernoulli cluster. The aim of this section is to introduce an adapted framework and to recall, without proofs, the results needed in this paper.

\subsection*{Grid structure of $\Z^d$} In the following, $d \geq 2$.
We denote by $\Z^d$ the 
graph whose set of vertices is 
$\Z^d$, and where we put a non-oriented edge between each pair $\{x,y\}$ of
\textit{neighbor} points in $\Z^d$, \ie points  whose Euclidean distance is equal to $1$. This set
of edges is denoted by $\Ed$. A \textit{(simple) path} is a sequence $\gamma=(x_1,x_2,\ldots,x_n,x_{n+1})$ of distinct points such that $x_i$ and $x_{i+1}$ are
neighbors and $e_i$ is the edge between $x_i$ and $x_{i+1}$. 
The number $n$ of edges in $\gamma$ is called the \textit{length} of
$\gamma$ and is denoted by $|\gamma|$. 

For any set $X$, and $u\in\Zd$, we define the \textit{translation operator}
$\theta_u$  
on $X^{\Ed}$ by the relation
$$\forall\omega\in X\quad\forall e\in \Ed\quad (\theta_u\omega)_e=\omega_{u.e},$$
where $u.e$ denotes the natural action of $\Zd$ on $\Ed$: if $e=\{a,b\}$, then
$u.e=\{a+u,b+u\}$. 

\subsection*{Assumptions and construction of  first-passage percolation}
$\;$
 
Denote by $p_c(d)$ the critical threshold for Bernoulli percolation on the
edges $\Ed$ of $\Z^d$, and choose $p\in (p_c,1]$.
On $\Omega_E=\{0,1\}^{\Ed}$, consider the measure $\P_p$:
$$
\mbox{on }\Omega_E=\{0,1\}^{\Ed}, \; \; 
\P_p=(p\delta_1+(1-p)\delta_0)^{\otimes \Ed}.
$$
A point $\omega$ in $\Omega_E$
is a  
\textit{random environment} for  first-passage percolation. An edge $e
\in \Ed$ is said to be 
\textit{open} in the environment $\omega$ if $\omega_e=1$, and \textit{closed}
otherwise. A path is said to be \textit{open} in the environment $\omega$
if all its edges are open in $\omega$. The \textit{clusters} of a
environment $\omega$ are the connected 
components of the graph induced on $\Z^d$ by the open edges in
$\omega$. As $p> p_c(d)$, there  almost surely exists one 
and only one infinite cluster, denoted by $C_\infty$.
On $\Omega_S=(\R_+)^{\Ed}$, consider a probability measure $\Snu$ such that:
$$
\mbox{on }\Omega_S=(\R_+)^{\Ed}, \; \; 
\Snu \mbox{ is stationary and
  ergodic}
$$
with respect to the previously introduced family of translations of the
grid. We suppose moreover that $\Snu$ satisfies integrability and dependence
conditions:
\begin{eqnarray}
& & m=\miniop{}{\sup}{e\in\Ed}\int \eta_{e}\
d\Snu(\eta)<+\infty. \label{sup-des-moments} \\
& & \exists \alpha>1, \; \exists A,B>0 \text{ such that }\forall
\Lambda \subseteq  \Ed, \; \Snu\left(\sum_{e\in\Lambda} \eta_i\ge
  B|\Lambda|\right)\le\frac{A}{| \Lambda|^{\alpha}}. \label{Halpha}
\end{eqnarray}
For instance, if $\Snu$ is the product measure $\nu^{\otimes\Ed}$,
assumption (\ref{Halpha}) follows from the Marcinkiewicz-Zygmund inequality
as soon as the passage time of an edge has a moment of order strictly
greater than $2$ -- see e.g. Theorem~3.7.8 in \cite{stout}.

Our probability space will then be $\Omega=\Omega_E\times \Omega_S$. 
A point in $\Omega$ will be
denoted $(\omega, \eta)$, with $\omega$ corresponding to the environment,
and $\eta$ assigning to each edge a non-negative \textit{passage time} which
represents the time needed to cross the edge. 
The final probability is:
$$
\mbox{on } \Omega=\Omega_E\times \Omega_S, \; \; \P=\P_p\otimes \Snu. 
$$

In the context of first-passage
percolation, as we are interested in asymptotic results concerning travel
time from the origin to points that tend to infinity, it is natural to
condition $\P_p$ on  the event that $0$ is in the
infinite cluster: 
$$\Pcond_p(.)=\P_p(.|0\in C_\infty) \; \; \mbox{ and } \; \;
\Pcond=\Pcond_p\otimes \Snu.$$ 
For $B\in\bor[\Omega_E]$, with $B\subset\{0\communique\infty\}$ and
$\P_p(B)>0$, we will also define the probability measure $\Pcond_B$ by 
$$\forall C\in\bor[\Omega]\quad \Pcond_B(C)=\frac{\P(C\cap
  (B\times\Omega_S))}{\P_p(B)}.$$ 

\subsection*{Examples}
The previous assumptions of the generalized first-passage percolation model
include:
\begin{itemize}
\item The case of classical \iid first-passage percolation: take $p=1$, \ie
all the 
edges of $\Z^d$ are open, and $\Snu=\nu^{\otimes \Ed}$, where $\nu$ is a
probability measure on $\R_+$.

\item The case of classical \iid first-passage percolation, but allowing the
passage times to take the value $\infty$ with positive probability:
take $p_c(d)<p<1$, a probability measure $\nu$ on $\R_+$, and set
$\Snu=\nu^{\otimes \Ed}$. This is equivalent to consider $p=1$ and
$\Snu=\tilde{\nu}^{\otimes \Ed}$, where $\tilde{\nu}$ is a 
probability measure on $\R_+ \cup \{\infty \}$ that charges $\infty$ with
probability $1-p$.

\item The case of stationary first-passage percolation, as considered by Boivin
in \cite{Boivin}: take $p=1$ and $\Snu$ a stationary probability measure.
\end{itemize}

\subsection*{The travel time}
The \textit{chemical distance} $D(x,y)$ between $x$ and
$y$ in $\Zd$ only depends on the Bernoulli percolation structure $\omega$
and is defined as follows: $D(x,y)(\omega)= \inf_{\gamma}|\gamma|,$
where the infimum is taken on the set of paths whose extremities are $x$
and $y$ and that are open in the environment $\omega$. 
It is of course only defined when $x$ and $y$ are in the same percolation
cluster, and represents then the minimal number of
open edges needed to link $x$ and $y$ in the environment
$\omega$. Otherwise, we set by convention $D(x,y)=+\infty$. 

For $(\omega,\eta)\in\Omega$, and $(x,y)\in\Zd\times\Zd$, we define the
\textit{travel time} from $x$ to $y$: 
$$d(x,y)(\omega,\eta)=\inf_{\gamma}d(\gamma)=\inf_{\gamma}\sum_{e\in
  \gamma}\eta_e,$$ 
where the infimum is taken on the set of paths whose extremities are $x$
and $y$ and that are open in the environment $\omega$. Of course
$d(x,y)=+\infty$ if and only  
if $D(x,y)=+\infty$.

A path $\gamma$ from $x$ to $y$ which realizes the distance $d(x,y)$
is called a \textit{finite geodesic}.
An infinite path $\gamma=(x_i)_{i\ge 0}$ is called a  \textit{semi-infinite geodesic} if $(x_n,x_{n+1},\dots,x_p)$ is a finite geodesic for every $n\le p$.

\subsection*{Directional asymptotic speed results}
In  classical first-passage percolation, we study, for each
$u\in\Zd\backslash\{0\}$, the travel time 
$d(0,nu)$ as $n$ goes to infinity. Here, as all points in $\Z^d$ are not
necessarily accessible from $0$, we must introduce the following definitions:
\begin{defi}
For each $u\in\Zd\backslash\{0\}$ and $B\in\bor[\Omega_E]$, let 
$$T^B_u(\omega)=\inf\{n\ge 1; \theta_{nu}\omega\in B\},$$ 
define the
associated random translation operator on $\Omega=\Omega_E \times \Omega_S$
$$
\Theta^B_u(\omega,\eta)
=\left( \theta_u^{T^B_u(\omega)}(\omega),\theta_u^{T^B_u(\omega)}(\eta)
\right)
$$
and the composed version 
$ \displaystyle
T^B_{n,u}(\omega)=\sum_{k=0}^{n-1} T^B_u \left( (\Theta_u^B)^k\omega \right).
$
\end{defi}
Note that $T^B_u$ only depends on the environment  $\omega$,
and not on the 
passage times $\eta$, whereas the operator $\Theta^B_u$ acts on the whole
configuration $(\omega,\eta)$. The next step is to study the asymptotic
behavior of such quantities:

\begin{lemme} \label{lemme_invariance}
$\Theta^B_u$ is a $\Pcond$-preserving transformation, 
is ergodic for $\Pcond$ and
$$\F T_u^B=\frac1{\P_p(B)} \; \; \mbox{ and } \; \;
\frac{T^B_{n,u}}n\to\frac1{\P_p(B)}\quad\Pcond\as$$ 
\end{lemme}

\begin{proof}
The idea is to prove that classical ergodic theorems can be applied.
\end{proof}

We turn now to the study of the quantity analogous to $d(0,nu)$ in the
classical first-passage percolation:
\begin{lemme}
\label{lemme-mu}
Let $B\in\bor[\Omega_E]$, with $B\subset\{0\communique\infty\}$ and
$\P_p(B)>0$. For $u\in\Zd\backslash\{0\}$, there exists a constant
$f^B_u\ge 0$ such that 
$$
\frac{d \left( 0,T^B_{n,u}(\omega)u \right) (\omega,\eta)}n \to
f^B_u\quad\Pcond_B\as$$ 
The convergence also holds in $L^1(\Pcond_B)$.
Moreover, $f^B_u\le \E_{\Pcond_B} d(0,T^B_u u)<+\infty$.
\end{lemme}

\begin{proof}
These results are proved with full details when  $B=\{0\communique\infty\}$
in \cite{garet-marchand}. Since the proof
is essentially the same, we omit it.  
\end{proof}

Now, for each $u\in\Zd\backslash\{0\}$, we define the asymptotic speed
in the direction $u$ by 
 $$\mu(u)=\P_p(0\communique\infty)f^A_u$$ for the choice 
 $A=\{0\communique\infty\}$. We also define $\mu(0)=0$. 

\begin{coro}
\label{coro-mu}
Let $B\in\bor[\Omega_E]$, with $B\subset\{0\communique\infty\}$ and
$\P_p(B)>0$. For $u\in\Zd\backslash\{0\}$, we have:
\begin{eqnarray*}
\frac{d(0,(T^B_{n,u}(\omega)u)(\omega,\eta)}n
& \to & \frac{\mu(u)}{\P_p(B)}\quad\Pcond_B\as \\
\frac{d(0,(T^B_{n,u}(\omega)u)(\omega,\eta)}{T^B_{n,u}(\omega)}
& \to & {\mu(u)}\quad\Pcond_B\as
\end{eqnarray*}
\end{coro}

\begin{proof}
We use the fact that $\left( \frac{d(0,T^B_{n,u}u)}{T^B_{n,u}}
\right)_{n\ge 0}$, as a 
subsequence of  $\left( \frac{d(0,T^A_{n,u}u)}{T^A_{n,u}} \right)_{n\ge
  0}$, admits 
the same almost sure limit $\mu(x)$, and the lemma~\ref{lemme_invariance}.
\end{proof}

In \cite{garet-marchand}, it has been 
proved that $\mu$ enjoys the properties that are usual in classical \iid
first-passage percolation: $\mu$ is a semi-norm. In classical \iid
first-passage percolation with passage time law $\nu$, it is well-known
that $\mu$ is a norm as soon as  $\nu(0)<p_c(d)$. In
~\cite{garet-marchand}, we gave a long 
discussion about conditions on $\Snu$ implying that $\mu$ is a norm.
Particularly, if $\Snu$ is a product measure $\nu^{\otimes\Zd}$,
$\mu$ is a norm as soon as $p\nu(0)<p_c(d)$.


\section{Coexistence result}

Consider the first-passage percolation model on $\Z^d$  previously
introduced. For every pair $x$ and $y$ of distinct points in $\Zd$, say that 
the event $\Coex(x,y)$ happens if  
$$\{z\in \Z^d; d(x,z)<d(y,z)\}\text{ and } \{z\in \Z^d; d(x,z)>d(y,z)\}$$
are both infinite sets. 

The goal of the paper is to prove that for every pair of distinct points $x,y\in\Zd$, $\P(\Coex(x,y))>0$. Our proofs always require the assumption that $\mu$ is not identically null and we guess that this assumption is close to be optimal. Let us detail a particular case where $\mu=0$ and coexistence never occurs.
Suppose that $d=2$ and $\Snu=\nu^{\otimes\Ed}$, with $p\nu(0)>p_c(2)=\frac12$. In this case, $\mu$ is identically null, as previously noted.
Consider two distinct points $x,y\in\Z^2$. Since $p\nu(0)>p_c(2)$, there almost surely exists an infinite cluster of open edges with passage time zero. It is known that in dimension 2, the supercritical infinite cluster almost surely contains a circuit that surrounds $x$ and $y$ and disconnects them from infinity -- see Harris~\cite{Harris} or for instance Grimmett's book~\cite{grimmett}. Clearly, the points in this circuit are equally $d$-distant from $x$ (\resp $y$). So, if $x$ reaches the circuit before $y$, it necessarily also reaches every point outside the circuit before $y$.  Similarly, if $x$ and $y$ reach the circuit 
at the same time, all the points outside the circuit will also be reached at the same time by $x$ and $y$.
In both cases, coexistence does not occur.

The next theorem gives conditions that ensures that 
coexistence possibly occurs for some (random) $x,y$.

\begin{theorem}
\label{yarfyarf}
Let $d \geq 2$, $p>p_c(d)$, $\Snu$ a
  stationary ergodic probability measure on $(\R_+)^{\Ed}$
  satisfying (\ref{sup-des-moments}) and (\ref{Halpha}), and $\mu$ be the related semi-norm describing the directional asymptotic speeds. 

Let $B\in\mathcal{B}(\Omega_E)$, with $B\subset\{0\communique\infty\}$ and
$\P_p(B)>0$, and $y\in\Zd$. We have:
$$
\mbox{if }\E d(0,T^B_{1,y}y)<\frac{2\mu(y)}{\P_p(B)}, \; 
\mbox{ then } \Pcond_B(\Coex(0,T^B_{1,y}y))>0.
$$
Moreover, if $x\in\Zd$ is such that $\mu(x)>0$, then $y=rx$ satisfies to
the previous condition provided that $r$ is large enough. 
\end{theorem}


Note that when $p=1$, which corresponds to classical first-passage percolation, we can take $B=\Omega_E$, and then $T^B_{1,y}y$ is simply equal to $y$.

\begin{figure}[h]
\label{fig-sym}
\setlength{\unitlength}{1cm}
\begin{picture}(8,1.8)
\put(0,0.3){\circle*{0.1}}\put(1,0.3){\circle*{0.1}}\put(2.5,1.3){\circle*{0.1}}
\put(-0.1,0.5){$A$}\put(0.9,0.5){$B$}\put(2.4,1.5){$M_n$}

\put(10,1.3){\circle*{0.1}}\put(9,1.3){\circle*{0.1}}\put(7.5,0.3){\circle*{0.1}}
\put(9.9,1.5){$J_n$}\put(8.9,1.5){$I_n$}\put(7.4,0.5){$O$}
\end{picture}
\caption{The symmetry argument}
\end{figure}
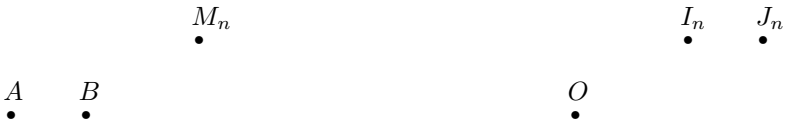

Before beginning the proof, we want to describe an elementary and  clever
trick used by  Pemantle and H{\"a}ggstr{\"o}m in~\cite{Haggstrom-Pemantle-1} that will
also be useful here. 
Consider figure~\ref{fig-sym}.
The left picture deals with our problem: if we prove that when
$M_n$ goes to the infinity on the right (\resp on the left), then $M_n$ is
infinitely often closer (\resp more distant) from $B$ than from $A$ with a
probability strictly larger than $0.5$, then coexistence holds with
positive probability. 

Now consider the right picture: for fixed $n$, 
$(d(0,J_n),d(0,I_n))$ has the same law that $(d(A,M_n),d(B,M_n))$, so 
the event $\{d(0,J_n)>d(0,I_n)\}$ occurs with the same probability as
the event $\{d(A,M_n)>d(B,M_n)\}$.
So, if we show that $\P(d(0,J_n)>d(0,I_n))>1/2$ for infinitely many  $n$,
the result is proved.  

As H{\"a}ggstr{\"o}m and Pemantle said, the idea is that there are sites
arbitrarily far away from the origin which strongly feel from which source the
infection is coming. Their \emph{modus operandi}, in the case of \iid exponentials on
$\Z^2$, was to control the infection rate ``from the right to the left'' and
the infection rate ``from the left to the right''. The main idea of the
proof which follows is that the advantage of the 
closest source can be quantified using the existence of a directional  asymptotic speed in first-passage percolation. Concretely,
we will use the law of $d(A,M_n)-d(B,M_n)$ (in fact, the law of
$d(0,J_n)-d(0,I_n)$) instead  of those of   $\{d(A,M_n)>d(B,M_n)\}$ (or
$\{d(0,J_n)>d(0,I_n)\}$) to carry the information. 

\begin{proof}(\textit{Theorem~\ref{yarfyarf}}) Choose $y \in \Z^d
  \backslash \{0\}$ such
  that
\begin{equation}
\label{undouble}
\E d(0,T^B_{1,y}y)<\frac{2\mu(y)}{\P_p(B)}. 
\end{equation}
Let us note
\begin{eqnarray*}
S_0 & = & \miniop{}{\limsup}{\|z\|_1\to
  +\infty}\{d(0,z)<d(T^B_{1,y}y,z)<+\infty\}, \\
S_1 & = & \miniop{}{\limsup}{\|z\|_1\to +\infty}\{ +\infty>
d(0,z)>d(T^B_{1,y}y,z)\}.
\end{eqnarray*}
It is obvious that $\Coex(0,T^B_{1,y}y)=S_0\cap S_1$. Intuitively, one
expects that the difference between $d(0,z)$ and $d(T^B_{1,y}y,z)$ will
be more important if $z \in \R y$, and we will effectively consider such $z$.
For the convenience of the reader, we also note, for $n\in\N$ and $x\in\Zd$,
$\tilde{T}_{n,x}=T^B_{n,x} x$. Define $\tilde{T}_{0,x}=0$, and for $n \geq 0$, 
$$
\begin{cases}
X_n=d(0,\tilde{T}_{n,y})-d( \tilde{T}_{1,y} ,\tilde{T}_{n,y}),\\ 
X'_n=d(\tilde{T}_{1,y} ,\tilde{T}_{n,-y})-d(0,\tilde{T}_{n,-y}).
\end{cases}
$$ 
By triangular inequality, one has 
$|X_n|\le d(0,\tilde{T}_{1,y})$ and $|X'_n|\le d(0,\tilde{T}_{1,y})$.

Note that for $\omega\notin S_0$, $X_n(\omega)\le 0$ as soon as $n$ is
large enough, whereas
for $\omega\notin S_1$, $X'_n(\omega)\le 0$ for large $n$.
It follows that for $\omega\notin S_0\cap S_1$, $X_n(\omega)+X'_{n-1}(\omega)\le
d(0,\tilde{T}_{1,y})(\omega)$ for large $n$. 
Let us define 
$$
Q_n=\sum_{k=1}^n (X_k+X'_{k-1}),\;\;  Z_n=\frac{Q_n}n\;\;\text{ and
  }\;\;Z=\miniop{}{\limsup}{n\to +\infty} Z_n.
$$
The previous remark implies easily that
\begin{equation}
\label{la_belle}
\forall\omega\notin S_0\cap S_1\quad Z(\omega)\le d(0,\tilde{T}_{1,y})  (\omega).
\end{equation}
By lemma~\ref{lemme-mu}, $d(0,\tilde{T}_{1,y})$ is integrable under $\Pcond_B$.
Since $|Z_n|\le  d(0,\tilde{T}_{1,y})$, it follows (for instance by Fatou's lemma) that 
$$\E_{\Pcond_B} Z=\E\miniop{}{\limsup}{n\to +\infty} Z_n\ge
\miniop{}{\limsup}{n\to +\infty} \E_{\Pcond_B} Z_n.
$$
Since $d( \tilde{T}_{1,y}
,\tilde{T}_{n,y})=d(0,\tilde{T}_{n-1,y})\circ\Theta^B_{y},$ it follows 
from the invariance of $\Pcond_B$ under $\Theta^B_{y}$ that
\begin{eqnarray*}
\E_{\Pcond_B} X_n
& = & \E\left( d(0,\tilde{T}_{n,y})-d( \tilde{T}_{1,y}
,\tilde{T}_{n,y})\right) \\
& = & \E_{\Pcond_B} d(0,\tilde{T}_{n,y})-\E_{\Pcond_B}
d(0,\tilde{T}_{n-1,y}).
\end{eqnarray*}
Then, it follows that $\E_{\Pcond_B}(X_1+X_2+\dots+X_n)=\E_{\Pcond_B} d(0,\tilde{T}_{n,y})$.
Similarly, as $d(\tilde{T}_{1,y}
,\tilde{T}_{n,-y})=d(0,\tilde{T}_{n+1,-y})\circ\Theta^B_{y},$ 
\begin{eqnarray*}
\E_{\Pcond_B} X'_n
& = & \left(\E_{\Pcond_B} d(\tilde{T}_{1,y}
,\tilde{T}_{n,-y})-d(0,\tilde{T}_{n,-y})\right) \\
& = & \E_{\Pcond_B} d(0,\tilde{T}_{n+1,-y})-\E_{\Pcond_B}
d(0,\tilde{T}_{n,-y}), 
\end{eqnarray*}
and $\E_{\Pcond_B}(X'_0+X'_1+\dots+X'_{n-1})=\E_{\Pcond_B} d(0,\tilde{T}_{n,-y})=\E_{\Pcond_B} d(0,\tilde{T}_{n,y}),$
using for the last equality the fact that $\Pcond_B$ is invariant under $(\Theta^B_{y})^n$ and the fact that a distance is symmetric.

Then, $\E_{\Pcond_B} Z_n=\frac{2\E_{\Pcond_B}
  d(0,\tilde{T}_{n,y})}n$. Since, via Corollary \ref{coro-mu},
$\frac{\E_{\Pcond_B} d(0,\tilde{T}_{n,y})}n$ converges to
$\frac{\mu(y)}{\P_p(B)}$, it follows that  
\begin{equation}
\label{minoration}
\E_{\Pcond_B} Z\ge \frac{2\mu(y)}{\P_p(B)}.
\end{equation}

Putting together (\ref{undouble}), (\ref{la_belle}) and
(\ref{minoration}), we see that 
$\Pcond_B(S_0\cap S_1)=0$ -- or equivalently $\Pcond_B((S_0\cap S_1)^c)=1$ --
would yield to a contradiction. This concludes the proof of the first
assertion.

The second assertion is a direct consequence of Corollary \ref{coro-mu}.
\end{proof}

One can be a bit perplexed by the fact that the position of the source which
may coexist with a source at the origin is a random variable.
The goal of the next result is to come back to deterministic
sources. Intuitively, one can guess that the larger the distance between
the two sources is, the higher the probability of coexistence will be. 
This is the spirit of the next result.

\begin{theorem}
\label{coex} 
Under the same assumptions as in Theorem \ref{yarfyarf}, 
suppose moreover that $\mu$ is not 
identically null. Then, we have:
\begin{itemize}
\item For $x\in\Zd$ with $\mu(x)\ne 0$, there is an infinite set of odd
  values for $n\in\N$\\ such that $\P (\Coex(0,nx))>0$. 
\item $\P(\exists x,y\in\Zd, \; \Coex(x,y))=1$.
\end{itemize}
\end{theorem}

Let us say a word on the unexpected apparition of odd integers. Of course,
the result would be the same with the set of integers and generally, this 
additional constraint does not bring much. 
Nevertheless, one will see later  that, in the competition context,  this additional property sometimes
prevents the two infections from reaching a point at the very same time; it
will also plays a fundamental role in the proof of Theorem~\ref{temps-entier}. 

\begin{proof}
Let $x\in\Zd$ be such that $\mu(x)>0$ and $N\in\N$.
Let $A=\{0\communique\infty\}$ and $B=A\cap \{T^A_{-x}\text{ is odd}\}$.
We have, from the FKG inequalities,
$$
\Pcond_p(B)\ge
\Pcond_p(T^A_{-x}=1)=\Pcond_p(-x\communique\infty)\ge
\P_p(-x\communique\infty)>0. 
$$ 
By lemma~\ref{lemme-mu} and Theorem~\ref{coro-mu}, $\frac{\E
  d(0,T^B_{1,rx})}{r}$ tends to 
$\frac{\mu(x)}{\P_p(B)}$, so we can find an odd integer $r\ge N$ with 
$\frac{\E d(0,T^B_{1,rx}rx)}{r}<\frac{2\mu(x)}{\P_p(B)}$.
By Theorem~\ref{yarfyarf}, one has $\Pcond(S_0\cap S_1)>0$.

By its definition,  $T^B_{1,rx}$ almost surely takes its values in the set of
non-negative odd integers. Then, we can  write
$$\Pcond_B(S_0\cap S_1)=\sum_{k\text{ odd}}\Pcond_B(S_0\cap
S_1\cap\{T^B_{1,rx}=k\}).$$ 
Then, there exists an odd integer $k\in\N$, with $\Pcond_B(S_0\cap
S_1\cap\{T^B_{1,rx}=k\})>0$. 
So, if we note $n=kr$, we have $n\ge r\ge  N$, $n$ is odd  and
$$\P(\Coex(0,nx))\ge \P_p(B)\Pcond_B(S_0\cap S_1\cap\{T^B_{1,rx}=k\})>0.$$

The second point is a consequence of the ergodicity assumption.
\end{proof}

  
\section{Mutual
unbounded growth and existence of two disjoint geodesics for
diffuse passage times}


 The aim of this section is to
prove the possibility of coexistence in general two-type
first-passage percolation, and 
to study the existence  of two semi-infinite geodesics in the corresponding
spanning tree.

We first give our assumptions on the law of the passage times $\Snu$ and
define the two-type competition model.

\subsection*{Assumptions.} We consider first-passage percolation on $\Z^d$,
with $d \geq 2$. The open edges are given by a
realization of a Bernoulli percolation on the edges $\Ed$ of $\Z^d$ with
parameter $p\in (p_c(d),1]$:
$$
\mbox{on }\Omega_E=\{0,1\}^{\Ed}, \; \; 
\P_p=(p\delta_1+(1-p)\delta_0)^{\otimes \Ed}.
$$
The passage times of the edges are given by a probability measure $\Snu$:
$$
\mbox{on }\Omega_S=(\R_+)^{\Ed}, \; \; 
\Snu \mbox{ is stationary and
  ergodic.}
$$
Finally, we consider the product measure $\P=\P_p \otimes \Snu$ on $\Omega_E
\times \Omega_S$.
We also need two distinct initial sources $s_1,s_2$ in $\Z^d$.

\begin{lemme}
\label{lem_defi_geod}
If $x \in \Z^d$ is such that
$d(s_1,x)<d(s_2,x)$ then for every $y$ in an optimal path realizing
$d(s_1,x)$, we have $d(s_1,y)<d(s_2,y)$.

\end{lemme}

\begin{proof}
Denote by $\gamma(s_1,x)$ an
optimal path from $s_1$ to $x$, and suppose that
there exists $y \in \gamma(s_1,x)$ such
that $d(s_2,y)\leq d(s_1,y)$. Then, by triangular inequality,
$$d(s_2,x) \leq d(s_2,y) +d(y,x)\leq d(s_1,y) +d(y,x)$$
but as $y \in \gamma(s_1,x)$, $d(s_1,x)=d(s_1,y)+d(y,x)$ and then
$d(s_2,x) \leq d(s_1,x)$,
which is a contradiction.
\end{proof}

This allows us to define the following two-type first-passage percolation
model.

\begin{defi}
\label{defi_Ntype}
Under the previous assumptions, we set:
$$
A_1=\{x \in \Z^d,  d(s_1,x)<d(s_2,x)\}, \; \; \mbox{ and } \; 
A_2=\{x \in \Z^d,  d(s_2,x)<d(s_1,x)\}.
$$
$A_i$ is the set of sites in
$\Z^d$ that 
are finally infected by type $i$ infection.
The time of infection of $x \in \Z^d$ is $t(x)=\inf\{ d(s_i,x), 1 \leq i
\leq 2\}$. We say  
that $x$ is finally infected if $t(x)<\infty$. 

Note that the set of finally
infected points could be larger than the union of $A_1$ and $A_2$: we can not \emph{ a priori} exclude that a point $x$ could be reached simultaneously by the two infections, in which case we call it an infected point without defining an  infection type. 

We say that the two infections \textit{mutually grow unboundedly} if the
two  sets $A_1$ and $A_2$ are both infinite.
\end{defi}

Thanks to lemma \ref{lem_defi_geod}, the
mutual unbounded growth 
of a two-type first-passage percolation starting from $s_1,s_2$ is equal to
the event $\Coex(s_1,s_2)$ defined in section~3. 
Note also that $A_1$ and $A_2$ are connected sets, and if they are
infinite, by  
a classical compactness argument, one
can find,  from each $x \in A_i$, a semi-infinite geodesic which is completely in  $A_i$. 

\vspace{0.2cm}
We add in this section assumptions that will ensure the uniqueness of 
optimal 
paths in  first-passage percolation under $\P_p \otimes \Snu$. If
$\Lambda$ is a finite subset of $\Ed$, denote by $\mathcal F_{\Lambda^c}$
the $\sigma$-algebra generated by $\{(\omega_e, \eta_e), e \notin
\Lambda\}$. We
suppose 
\begin{equation}
\label{diffus}
\forall \Lambda \mbox{ finite subset of } \Ed, \forall e \in \Lambda,
\forall a \in \R_+, \; \; \Snu(\eta_e=a|\mathcal F_{\Lambda^c})=0.
\end{equation}

\begin{lemme}
\label{lem-diffus}
Under the additional assumption (\ref{diffus}), we have:
\begin{enumerate}
\item
If $\gamma$ and $\gamma'$ are paths that differ at least from one edge,
then
$$\P \left( d(\gamma)=\sum_{e \in \gamma} \eta_e= d(\gamma')=\sum_{e \in \gamma'} \eta_e< \infty
  \right) =0.$$
\item
For every $x$ and $y$ in $\Z^d$, the travel time $d(x,y)$, when it is
finite,  is a.s. realized by a
unique path.
\item
For every $\alpha \in \R$, if $x,x',y,y'$ are distinct points in $\Z^d$,\\ $\P(d(x,y)-d(x',y')=\alpha)=0.$
\item
If $\displaystyle\inf_{1 \leq i \leq 2} d(s_i,x)$ is finite, it is realized by a unique
source $s_i$.
\end{enumerate}
\end{lemme}

\begin{proof}
These are classical and not too difficult consequences of
assumption 
(\ref{diffus}). 
\end{proof}

Now,
if $t(x)<\infty$, then $x$ is reached first by a unique infection; the
path of infection $\gamma(x)$ is the unique path from the corresponding
source to $x$ that realizes $t(x)$. The set of eventually infected
points is in this case the union of $A_1$ and $A_2$. In other words, we can
define uniquely, for each eventually infected point, its type of infection
and its optimal path.

The union of $(\gamma(x))_{x \in \Z^d, t(x)<\infty}$ is then a random forest of two
trees $T(s_1)$ and $T(s_2)$,  respectively rooted at $s_1$ and $s_2$ and
 respectively spanning $A_1$ and $A_2$. 
 
The next theorem ensures the irrelevance of the positions of the two initial sources in
determining whether mutual unbounded growth occurs with positive
probability or not. Its proof is based on a modification of the
configuration around the sources, sufficiently strong to change the initial
sources, and sufficiently slight to ensure that some geodesics are not
modified outside a finite box.

\begin{lemme}
\label{sources}
Consider $\Z^d$, with $d \geq 2$ and $p \in (p_c(d),1]$. Choose a stationary ergodic probability measure $\Snu$
on $\Omega_S=(\R_+)^{\Ed}$ satisfying to the non-atomic  assumption (\ref{diffus}) and to:
\begin{equation}
\label{energiefinie}
\forall \Lambda \mbox{ finite subset of } \Ed, \forall e \in \Lambda,
\forall \varepsilon>0, \; \; \Snu(\eta_e\leq \varepsilon |\mathcal F_{\Lambda^c})>0 \; a.s.
\end{equation}
If $p=1$, we add the assumption that the support of the passage time is
 conditionally unbounded: 
\begin{equation}
\label{suppnonborne}
\forall \Lambda \mbox{ finite subset of } \Ed, \forall e \in \Lambda,
\forall M>0, \; \; \Snu(\eta_e\geq M |\mathcal F_{\Lambda^c})>0 \; a.s.
\end{equation}
Then if $s_1,s_2$ and $s'_1,s'_2$ are two pairs of distinct points in $\Z^d$,
$$\P(\Coex(s_1,s_2))>0 \Leftrightarrow
\P(\Coex(s'_1,s'_2))>0.$$ 
\end{lemme}

Let us comment the two  assumptions (\ref{energiefinie}) and
(\ref{suppnonborne}). 
They have the form of finite energy properties, which are usual in modification
arguments: it enables to force the occurrence of a wished event inside a finite box. But they also enable the passage time of an edge to take as small -- and as large when $p=1$ -- values as we like.
This is rather a technical assumption that could probably be relaxed. 
For instance, assumptions   (\ref{energiefinie}) and
(\ref{suppnonborne}) are satisfied for $\Snu=\nu^{\otimes\Ed}$ with $\nu$ is 
equivalent to Lebesgue's measure on $\R_{+}$.

The next result says that the mutual unbounded growth in the
two-type first-passage percolation model and the existence of
two distinct semi-infinite geodesics in the embedded   spanning tree in the
corresponding first-passage 
percolation model are equivalent.
 
\begin{lemme}
\label{equivalence-coex-geod}
Under the same assumptions as in Theorem \ref{sources}, 
$$
\begin{array}{ll}
& \exists s_1,s_2 \in \Z^d \mbox{ such that }\P(\Coex(s_1,s_2))>0 \\
\Leftrightarrow & \P \left(
  \begin{array}{c}
    \mbox{there exist two edge-disjoint semi-infinite geodesics}\\
    \mbox{in the infection tree rooted in }0
  \end{array}
\right)>0.
\end{array}$$
\end{lemme}

Combining these results with Theorem~\ref{coex}, we obtain:
\begin{theorem}
\label{geod}
Consider $\Z^d$, with $d \geq 2$ and $p \in (p_c(d),1]$. Choose a stationary ergodic probability measure $\Snu$
on $\Omega_S=(\R_+)^{\Ed}$ satisfying to the
integrability assumptions (\ref{sup-des-moments}), (\ref{Halpha}) and 
such that the related semi-norm $\mu$ describing the directional asymptotic speeds is not identically null.

Suppose  moreover that $\Snu$
satisfies to 
the 
non-atomic  assumption
$$\forall \Lambda \mbox{ finite subset of } \Ed, \forall e \in \Lambda,
\forall a \in \R_+, \; \; \Snu(\eta_e=a|\mathcal F_{\Lambda^c})=0,$$
 and to the 
 finite energy property
$$\forall \Lambda \mbox{ finite subset of } \Ed, \forall e \in \Lambda,
\forall \varepsilon>0, \; \; \Snu(\eta_e\leq \varepsilon |\mathcal F_{\Lambda^c})>0 \; a.s.
$$
If $p=1$, we add the assumption that the support of the passage time is
 conditionally unbounded: 
\begin{equation*}
\forall \Lambda \mbox{ finite subset of } \Ed, \forall e \in \Lambda,
\forall M>0, \; \; \Snu(\eta_e\geq M |\mathcal F_{\Lambda^c})>0 \; a.s.
\end{equation*}
\noindent Then,

1. $\displaystyle \forall x \neq y \in \Z^d, \; \P(\Coex(x,y))  >  0.$

\vspace{0.2cm}
2. $\displaystyle \P \left(
  \begin{array}{c}
    \mbox{there exist two edge-disjoint semi-infinite geodesics}\\
    \mbox{in the infection tree rooted in }0
  \end{array}
\right)  >  0.$
\end{theorem}

\subsection*{Examples}
1. Consider first-passage percolation on $\Z^d$, $d \geq 2$, with a family
$(t(e))_{e \in \Ed}$ of 
\iid non negative random variables with a non-atomic unbounded support containing $0$, for
instance an exponential law as in Richardson's model. 
Then 
\begin{itemize}
\item
For the two-type competition model, the probability of mutual
unbounded growth is positive for every pair of distinct sources in $\Z^d$.
\item
For the first-passage percolation model with one
source, the probability that the embedded spanning tree of $\Z^d$ has two
edge-disjoint infinite branches is positive.
\end{itemize}
These results were proved by H{\"a}ggstr{\"o}m and Pemantle in the pioneer work
\cite{Haggstrom-Pemantle-1} 
for Richardson's model in dimension 2. 
Our results positively answer   to the questions asked by H{\"a}ggstr{\"o}m and Pemantle about extensions of their coexistence result to higher dimensions and more general distributions for passage times.

2. Consider first-passage percolation on $\Zd$, $d \geq 2$, with a family
$(t(e))_{e \in \Ed}$ of 
\iid non negative random variables 
whose law has 
no atom excepted in $\infty$ (\ie edges can be closed with positive
probability)
and has $0$ in its support, for instance $t \sim p .\mathcal U_{[0,1]} +
(1-p).\delta_\infty$, with $p_c(d)<p<1$. Then
\begin{itemize}
\item
For the two-type competition model, the probability of mutual
unbounded growth is positive for every pair of distinct sources in $\Zd$.
\item
For the first-passage percolation model with one
source, the spanning tree of the infinite open cluster has two
edge-disjoint infinite branches with positive probability.
\end{itemize}

\subsection*{Remarks} We evoke  here some possible extensions of these
results. 

1. In the spirit of Deijfen and H{\"a}ggstr{\"o}m's work~\cite{Deijfen-Haggstrom}, we could have considered
competition models with \textit{fertile finite sets} as initial sources
rather than 
\textit{points}. As the argument is a \textit{local} modification argument
around the sources, our proof can be adapted to generalize the
irrelevance of the initial sources result:
if $S_1,S_2$ and $S'_1,S'_2$ are two pairs of fertile finite sets in $\Z^d$,
$$\P(\Coex(S_1,S_2))>0 \Leftrightarrow
\P(\Coex(S'_1,S'_2))>0.$$ 

2. Let us say a word on \textit{multitype} first-passage percolation. The
definitions concerning the two-type first-passage percolation can be
generalized in the obvious manner to consider a competition model between
$N$ infections starting from $N$ sources $s_1,s_2,...,s_N$ and trying to
invade the sites of $\Z^d$. In this context, the event
$\Coex(s_1,s_2,...,s_N)$ is defined as the event that there finally exist
an infinite set of infected points of each type of infection. Theorems
\ref{sources} and \ref{equivalence-coex-geod} can be proved in the 
same manner for $N$-type first-passage percolation. The only difficulty is to
ensure that the considered initial sources $s_1,s_2,...,s_N$ are
susceptible to give rise to a coexistence configuration: this means initial
sources $s_1,s_2,...,s_N$ for which it is possible to
find a family of $N$ infinite paths $(\Gamma_i)_{1 \leq i \leq N}$ such
that for every $i$, $\Gamma_i$
starts from $s_i$ and such that any two of these paths have no point in
common.  

Unfortunately, the coexistence result Theorem~\ref{geod} is not
available for $N$ sources, as it relies on Theorem~\ref{yarfyarf}, which is
only valid for two sources, and whose proof doesn't seem to be easy to
adapt to more sources.

\vspace{0.2cm}
We can now begin the proofs of these results. As the arguments of lemmas
\ref{sources} and \ref{equivalence-coex-geod} are very similar, we give the
proof of lemma~\ref{sources} 
in full details, and give only indications to adapt the
proof for the geodesics problem.

\begin{proof} \textit{(lemma \ref{sources}).}
Choose $s_1,s_2$ and $s'_1,s'_2$ two pairs of distinct points in
$\Z^d$ and denote by $\Lambda$ an hypercubic box in $\Z^d$ large enough
to contain 
$s_1,s_2$ and $s'_1,s'_2$. We also define
$\partial \Lambda  =  \{ x \in \Lambda, \; \exists y\notin \Lambda, \; \|x-y\|_1=1\}$.

By enlarging $\Lambda$ if necessary, we can assume that $s_1,s_2,s'_1,s'_2$ are
at a distance at least $3$ from $\partial \Lambda$. 
For an edge $e \in \Ed$, we say that
$e \in \Lambda$ if and only if  $\mbox{its two extremities are in }\Lambda
\mbox{ and at least one is not in } \partial \Lambda.
$
For a point $(\omega,\eta)$  in $\Omega=\Omega_E \times
\Omega_S=\{0,1\}^{\Ed} \times (\R_+)^{\Ed}$, 
$$
(\omega_\Lambda,\eta_\Lambda)=\{(\omega_e,\eta_e), e \in \Lambda\}
\; \; \mbox{ and }\; \;  
(\omega_{\Lambda^c}, \eta_{\Lambda^c})=\{(\omega_e, \eta_e), e \in \Ed
\backslash \Lambda\}. 
$$
For two point $x,y$ that are  in $\Lambda^{c}\cup \partial \Lambda$, we
define 
$d_{\Lambda^c}(x,y)(\omega)$ as the infimum, among all the paths $\gamma$
from $x$ to $y$   
\textit{ whose edges are not in }$\Lambda$, of $\sum_{e \in
  \gamma} \eta_e$, and, when this quantity is finite,
$\gamma_{\Lambda^c}(x,y)$ is the only path that realizes this infimum.

Suppose that $\P(\Coex(s_1,s_2))>0$. This is equivalent to say that
$$
\P \left(
\begin{array}{l}
  T(s_1) \mbox{ contains a infinite branch }\Gamma_1 \mbox{ starting from
    }s_1, \\ 
  T(s_2) \mbox{ contains a infinite branch }\Gamma_2 \mbox{ starting from }s_2.
\end{array}
\right)>0.
$$
Remember that the box $\Lambda$ has been chosen large enough to contain
$s_1,s_2$, so there exists on $\Gamma_1$ (\resp $\Gamma_2$) a last point $r_1$ (\resp $r_2$) to be in
$\partial \Lambda$. As $\partial \Lambda$ is finite, there must exist two
distinct 
points $r_1,r_2 \in \partial \Lambda$ such that:
\begin{equation}
\label{eq_branches}
\P \left(
\begin{array}{c}
  T(s_1) \mbox{ contains a infinite branch }\Gamma_1 \mbox{ starting from
    }s_1 \\ 
  \mbox{ and whose last point in } \partial \Lambda \mbox{ is }r_1, \\
  T(s_2) \mbox{ contains a infinite branch }\Gamma_2 \mbox{ starting from
    }s_2 \\
  \mbox{ and whose last point in } \partial \Lambda \mbox{ is }r_2.
\end{array}
\right)>0.
\end{equation}
Now, we must introduce the following events:
\begin{eqnarray*}
C_1 & = & \left\{
\begin{array}{c}
  \mbox{There is a simple path } (x_i^1)_{i \geq 1} \mbox{ in }\Lambda^c
  \mbox{ such that } \Vert x_1^1-r_1\Vert_1=1\\ \mbox{ and }
  \forall i \geq 1,
  \gamma_{\Lambda^c}(r_1,x_i^1)=(r_1,x_1^1,...,x_{i-1}^1 ,x_i^1).
\end{array}
\right\},\\
C_2 & = & \left\{
\begin{array}{c}
  \mbox{There is a simple path } (x_j^2)_{j \geq 1} \mbox{ in }\Lambda^c
  \mbox{ such that } \Vert x_1^2-r_2\Vert_1=1\\ \mbox{ and }
  \forall j \geq 1,
  \gamma_{\Lambda^c}(r_2,x_j^2)=(r_2,x_1^2,...,x_{j-1}^2 ,x_j^2).
\end{array}
\right\},\\
A_i^1 & = & \left\{ d_{\Lambda^c}(x_i^1,r_1)+d(r_1,s_1)<
d_{\Lambda^c}(x_i^1,r_2)+d(r_2,s_2)\right\},\\
A_j^2 & = & \left\{ d_{\Lambda^c}(x_j^2,r_2)+d(r_2,s_2)<
d_{\Lambda^c}(x_j^2,r_1)+d(r_1,s_1)\right\}. \\
\end{eqnarray*}

\begin{figure}
\includegraphics[height=5.5cm]{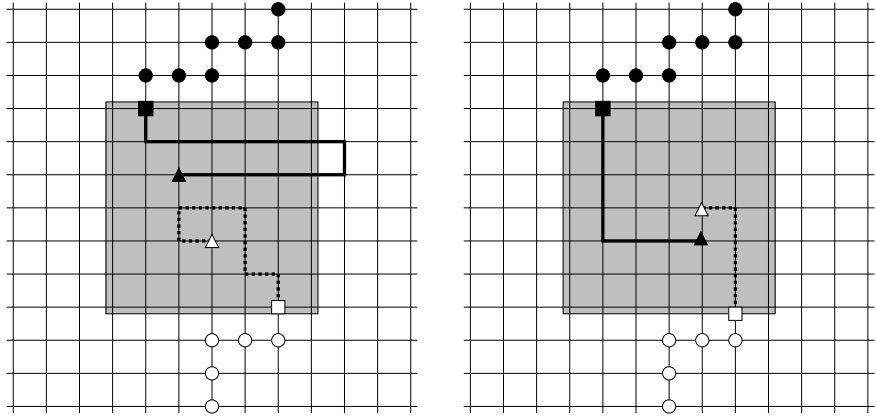}
\caption{Modification of the infection trees. On the left, competition with
  two sources $s_1$ (black triangle) and $s_2$ (white triangle): the box
  $B$ is in grey, the squares are the exiting points of the branches from
  $\Lambda$ (black for $r_1$ and 
  white for $r_2$) and the circles are the visible portions of the
  infinite branches outside $\Lambda$ (black for $x_1^1,x_2^1,...,x_7^1$ and
  white for 
  $x_1^2,...,x_5^2$). On the right, the configuration outside $\Lambda$ has not
  changed, but we forced the branches to follow $\gamma'_1$ (in black) and
  $\gamma'_2$ (dashed) and thus we changed the sources into
   $s'_1$, $s'_2$.}
\end{figure}

Let us prove that (\ref{eq_branches}) implies:
\begin{equation}
\label{eq_xy}
\P \left( C_1 \cap C_2 \cap \bigcap_{i \geq 1} A_i^1 \cap \bigcap_{j \geq 1}
  A_j^2 \right) >0.
\end{equation}
Indeed, suppose that the event in (\ref{eq_branches}) is realized. Then the
portion $(x_i^1)_{i \geq 1}$ of $\Gamma_1$ after $r_1$ is a good candidate
for $C_1$. As $\Gamma_1$
is a branch of the infection tree $T(s_1)$, one has $\forall i \geq 1, \gamma(r_1,x_i^1)=(r_1,x_1^1,...,x_{i-1}^1,x_i^1)$, which is stronger than the assertion needed for $C_1$.
Next, as $r_1$ is in $\Gamma_1$, for every $i \geq 1$ we know that $r_1$ is
a point of $\gamma(s_1, x_i^1)$, and the fact that $x_i^1$ is in $T(s_1)$
implies that 
$$d(s_1,r_1)+d(r_1,x_i^1)=d(s_1,x_i^1)< d(s_2,x_i^1)\leq d(s_2,r_2)+d(r_2,x^1_i).$$
As $r_1$ (\resp $r_2$) is the last point of $\Gamma_1$ (\resp $\Gamma_2$)
to be in $\Lambda$, we have $d(r_1,x_i^1)=d_{\Lambda^c}(r_1,x_i^1)$ (\resp 
$d(r_2,x_i^1)=d_{\Lambda^c}(r_2,x_i^1)$), and thus $A_i^1$ is realized. 
Doing the same for $C_2$ and $A_j^2$, we see that the event that appears in
(\ref{eq_branches}) is included in $C_1 \cap C_2 \cap \bigcap_{i \geq 1}
A_i^1 \cap \bigcap_{j \geq 1} A_j^2$, and thus (\ref{eq_xy}) is proved.

Now, as $C_1$ and $C_2$ are in $\mathcal F_{\Lambda^c}$, conditioning on
$\mathcal F_{\Lambda^c}$ gives: 
\begin{eqnarray*}
&  & \P \left( C_1 \cap C_2 \cap \bigcap_{i \geq 1} A_i^1 \cap \bigcap_{j
    \geq 1} A_j^2 \right) \\
& = & \int d\P(\omega_{\Lambda^c}, \eta_{\Lambda^c})
. \mathbf{1}_{C_1}(\omega_{\Lambda^c},\eta_{\Lambda^c} ) 
  . \mathbf{1}_{C_2}(\omega_{\Lambda^c},\eta_{\Lambda^c}) . \P \left(
    \left. (\omega, 
      \eta)  \in \bigcap_{i
  \geq 1} A_i^1 \cap \bigcap_{j \geq 1} A_j^2 \right| \mathcal
F_{\Lambda^c} \right), \\ 
\end{eqnarray*}
where we can also write:
\begin{eqnarray}
& & \P \left( \left. (\omega,\eta) \in \bigcap_{i
  \geq 1} A_i^1 \cap \bigcap_{j \geq 1} A_j^2 \right| \mathcal F_{\Lambda^c}
\right) \label{eq_inegalites}\\ 
& = & \P \left( \left. 
\begin{array}{l}
\forall i \geq 1, \;
d_{\Lambda^c}(x_i^1,r_1)-d_{\Lambda^c}(x_i^1,r_2)<d(r_2,s_2)-d(r_1,s_1) \\
\forall j \geq 1, \; 
d(r_2,s_2)-d(r_1,s_1)<d_{\Lambda^c}(x_j^2,r_1)-d_{\Lambda^c}(x_j^2,r_2)
\end{array}
\right| \mathcal F_{\Lambda^c} \right) \; a.s. \nonumber
\end{eqnarray}
Define 
\begin{eqnarray*}
m_1=m_1(\omega_{\Lambda^c},\eta_{\Lambda^c})& = & \sup_{i \geq 1}
(d_{\Lambda^c}(x_i^1,r_1)-d_{\Lambda^c}(x_i^1,r_2)) \\
m_2=m_2(\omega_{\Lambda^c},
\eta_{\Lambda^c}) & = & \inf_{j \geq 1}
(d_{\Lambda^c}(x_j^2,r_1)-d_{\Lambda^c}(x_j^2,r_2)). 
\end{eqnarray*} 
Now, we have 
\begin{equation}
\label{eq_m1m2}
\P \left( C_1 \cap C_2 \cap \bigcap_{i \geq 1} A_i^1 \cap \bigcap_{j
    \geq 1} A_j^2 \cap \{m_1<m_2\}
\right)>0.
\end{equation}
Indeed, thanks to (\ref{eq_inegalites})
and to lemma \ref{lem-diffus}.iii,
$$\P \left( \left. \bigcap_{i \geq 1} A_i^1 \cap \bigcap_{j \geq 1} A_j^2 
\cap \{ m_1 \geq m_2 \} \right| \mathcal F_{\Lambda^c}
\right)= \P \left( \left. d(r_2,s_2)-d(r_1,s_1) = m_1 \right| \mathcal
  F_{\Lambda^c} \right) =0 \; a.s. $$
and then the probabilities in (\ref{eq_m1m2}) and in (\ref{eq_xy}) are equal.
It is
also easy to see, if $m_1(\omega_{\Lambda^c},
\eta_{\Lambda^c})<m_2(\omega_{\Lambda^c}, \eta_{\Lambda^c})$, that we can
find $a_1,a_2,b_1,b_2 \in \R_+$ 
such that 
$$
b_1m_1<b_2m_2,\quad  a_2m_2<a_1m_1,\quad b_1-a_1=1,\quad\mbox{ and } b_2-a_2=1.
$$
Define also
\begin{eqnarray*}
M=M(\omega_{\Lambda^c}, \eta_{\Lambda^c})& = & \max \left\{ a_1m_1, b_2m_2\right\} \\
&+&\max_{\substack{z \in \partial\Lambda\\d_{\Lambda^c}(r_1,z)< \infty}} \left\{ d_{\Lambda^c}(r_1,z)
\right\} 
+\max_{\substack{z \in \partial \Lambda\\ d_{\Lambda^c}(r_2,z)< \infty}} \left\{ d_{\Lambda^c}(r_2,z)
\right\}. 
\end{eqnarray*}
Now, we build a set $G=G(\omega_{\Lambda^c}, \eta_{\Lambda^c})$ of
\textit{good} 
configurations $(\omega_{\Lambda}, \eta_\Lambda)$ inside $\Lambda$,
depending on the configuration outside $\Lambda$. First, since $\Lambda$ has been chosen large enough, it is possible to 
draw with the edges in $\Lambda$, a path $\gamma'_1$ that links $s'_1$ to $r_1$
and a path $\gamma'_2$ that links $s'_2$ to $r_2$ such that $\gamma'_1$ and
$\gamma'_2$ have no vertex and no edge in common. Denote by $|\gamma'_1|$
(\resp $|\gamma'_2|$) the number of edges in $\gamma'_1$
(\resp $\gamma'_2$). 
We define now $G$ as the set of $(\omega_{\Lambda}, \eta_\Lambda)$ that satisfy the
following conditions: 
\begin{enumerate}
\item
  $\forall e \in \gamma'_1, \; \; \omega_e=1 \; \; \mbox{ and } \; \;
  a_2m_2/|\gamma'_1|<\eta_e< a_1m_1/|\gamma'_1|$, 
\item
  $\forall e \in \gamma'_2,  \; \; \omega_e=1 \; \; \mbox{ and } \; \;
  b_1m_1/|\gamma'_2|<\eta_e< b_2m_2/|\gamma'_2|$, 
\item
\begin{itemize}
\item
if $p<1$, then $\forall e \in \Lambda \backslash (\gamma'_1 \cup \gamma'_2), \;
\; \omega_e=0$. 
\item
if $p=1$, then $\forall e \in \Lambda \backslash (\gamma'_1 \cup \gamma'_2), \; 
\; \eta_e> M$.
\end{itemize}
\end{enumerate}

Under the finite energy assumptions (\ref{energiefinie}) and
(\ref{suppnonborne}), on the event $\{m_1<m_2\}$, we have 
$\P(G(\omega_{\Lambda^c}, \eta_{\Lambda^c})|\mathcal F_{\Lambda^c})>0\ \as$, so
(\ref{eq_m1m2}) 
implies  
\begin{eqnarray}
\int  \mathbf{1}_{C_1}
. \mathbf{1}_{C_2} . \mathbf{1}_{\{m_1<m_2\}} . \P \left( G |\mathcal F_{\Lambda^c} \right) d\P 
& = & \P (C_1 \cap C_2 \cap \{m_1<m_2\} \cap G)>0. \label{eq_intposi}
\end{eqnarray}
Let us prove now that on the event $C_1 \cap C_2 \cap \{m_1<m_2\} \cap G$,
each of the two 
infection trees $T(s'_1)$ and $T(s'_2)$ contains one infinite branch, or in
other words 
$$C_1 \cap C_2 \cap \{m_1<m_2\} \cap G \subset \Coex(s'_1,s'_2).$$ 
Suppose then 
that $(\omega,\eta) \in C_1 \cap C_2 \cap \{m_1<m_2\} \cap G$. We have, in
 the configuration $(\omega,\eta)$: 
\begin{itemize}
\item
$a_2m_2< d(\gamma'_1)<a_1m_1$ thanks to
condition i. in the definition of $G$.
\item
$b_1m_1<  d(\gamma'_2)<b_2m_2$ thanks to
condition ii. in the definition of $G$.
\item
Thus, by difference and by the choice of $a_1,b_1,a_2,b_2$, we have \\
$m_1= m_1(b_1-a_1)< d(\gamma'_2)-d(\gamma'_1)<
m_2(b_2-a_2)=m_2.$
\label{eq_tempsintOK}
\end{itemize}
Moreover, as soon as a path $\gamma$ from $s'_1$ to $r_1$ differs from
$\gamma'_1$ by at least one edge, it must use an edge $e$ in
$\Lambda \backslash (\gamma'_1 \cup \gamma'_2)$, and this edge is either
closed or such that $\eta_e>M$ thanks to condition iii. in the definition of
$G$. Thus $d(\gamma)>M\geq a_1m_1
>d(\gamma'_1)$. Then, $\gamma'_1$ is the optimal path from $s'_1$ to
$r_1$. On the other hand, every path $\gamma$ from $s'_2$ to $r_1$ has to
use an edge $e$ in $\Lambda \backslash (\gamma'_1\cup \gamma'_2)$, and then
by the same argument  
$d(\gamma) >M\geq a_1m_1 >d(\gamma'_1)$, and then $d(s'_2, r_1)>
d(s'_1,r_1)$. Consequently,  
$$
r_1 \in T(s'_1)\; \; \mbox{ and } \; \; \gamma'_1 \subset T(s'_1)\; \;
\mbox{ and } \; \; d(s'_1,r_1)=d(\gamma'_1). 
$$
In the same manner, 
$$
r_2 \in T(s'_2)\; \; \mbox{ and } \; \; \gamma'_2 \subset T(s'_2)\; \;
\mbox{ and 
  } \; \; d(s'_2,r_2)=d(\gamma'_2). 
$$
Let us prove now that the path $(x_i^1)_{i \geq 1}$, given by $C_1$,  is
an infinite branch of 
$T(s'_1)$. Let $i \geq 1$, and let $\gamma$ be a path from $s'_1$ to
$x_i^1$ that doesn't exit from $\Lambda$ in $r_1$ but in $z \neq r_1$. If $\gamma$ is
the minimal path from $s'_1$ to $x_i^1$, we must have
\begin{eqnarray*}
&& d(s'_1,x_i^1) = d(\gamma) =d(s'_1,z)+d_{\Lambda^c}(z,x_i^1)<
d(s'_1,r_1)+d_{\Lambda^c}(r_1,x_i^1) \\
\mbox{and then} && d(s'_1,z)<
d(r_1,s'_1)+d_{\Lambda^c}(r_1,x_i^1)-d_{\Lambda^c}(z,x_i^1) \leq 
d(r_1,s'_1)+d_{\Lambda^c}(r_1,z) .
\end{eqnarray*}
The last inequality is just the triangular inequality for $d_{\Lambda^c}$.
But the path
$\gamma$ must then contain at least one edge $e$ in $\Lambda \backslash (
\gamma'_1\cup \gamma'_2)$,  and so such that $\eta_e>M$ or 
$\omega_e=0$, and 
then by definition of $M$, we must have
$$
d(s'_1,z) >M \geq a_1m_1+d_{\Lambda^c}(r_1,z) >
d(r_1,s'_1)+d_{\Lambda^c}(r_1,z),$$
which  contradicts  the previous inequality. Thus, the portion of the
branch between $s'_1$ 
and $x_i^1$ that is in $\Lambda$ exits $\Lambda$ in $r_1$: it  is exactly
$\gamma'_1$. It remains to see that the part of the branch from $s'_1$ to
$x_i$ that is not in $\Lambda$ is exactly $(r_1,x_1,...,x_{i-1},x_i)$. But this
condition is always satisfied in the event $C_1$. Thus the path $(x_i^1)_{i \geq 1}$, given by $C_1$,  is
an infinite branch of 
$T(s'_1)$, and in the same manner, the path $(x_j^2)_{j \geq 1}$, given
by $C_2$,  is 
an infinite branch of 
$T(s'_2)$. We have thus proved the desired inclusion $C_1 \cap C_2 \cap \{m_1<m_2\} \cap G \subset
\Coex(s'_1,s'_2)$. 

Now, (\ref{eq_intposi}) ensures that $\P( \Coex(s'_1,s'_2))>0$, which ends
the proof.
\end{proof}

\begin{proof}(\textit{Theorem \ref{equivalence-coex-geod}}).
To prove the direct implication,
the proof is the same with 
$s'_1=s'_2=0$. The only difference  is to take the two paths
$\gamma'_1,\gamma'_2$ rooted both in $0$, with no other point in common and
with no edge in common.

The converse implication can also be proved by an analogous modification
argument. 
\end{proof}



\section{Mutual unbounded growth and existence of two distinct geodesics
  for integer passage times} 


In the previous section, the law of the passage time of an edge was supposed
to admit no atom, and thus the minimal paths were unique. However, as seen
in definition~\ref{defi_Ntype}, it is still possible to define a
two-type competition model. Choose two
distinct sources $s_1,s_2 \in \Z^d$. We define exactly as previously the sets:
$$A_1=\{x \in \Z^d,  d(s_1,x)< d(s_2,x)\}, \mbox{ and } A_2=\{x \in \Z^d,
d(s_2,x)< d(s_1,x)\}, $$
and say that coexistence occurs, event denoted by  $\Coex(s_1,s_2)$, if these two sets are
infinite. 

As in the preceding section, the geodesics  inside
$A_1$ and $A_2$  will be useful to obtain
configurations that allow to prove the irrelevance of the initial sources
in determining whether coexistence happens with positive probability or not.

\subsection*{Assumptions.} We consider first-passage percolation on $\Z^d$,
with $d \geq 2$. The open edges are given by a
realization of a Bernoulli percolation on the edges $\Ed$ of $\Z^d$ with
parameter $p \in (p_c(d),1]$:
$$
\mbox{on }\Omega_E=\{0,1\}^{\Ed}, \; \; 
\P_p=(p\delta_1+(1-p)\delta_0)^{\otimes \Ed}.
$$
The passage times of the edges are given by a probability measure $\Snu$:
$$
\mbox{on }\Omega_S=(\R_+)^{\Ed}, \; \; 
\Snu \mbox{ is stationary and
  ergodic.}
$$
Finally, we consider the product measure $\P=\P_p \otimes \Snu$ on $\Omega_E
\times \Omega_S$. We say that $\Snu$ satisfies condition $(H_\alpha)$ if 
\begin{equation*}
\label{halphalavraie}
(H_{\alpha})\quad \exists A,B>0\ \text{such that }\quad\forall
\Lambda\subseteq  \Ed, \; \Snu\left(\eta\in\Omega_S;\sum_{e\in\Lambda} \eta_i\ge
  B|\Lambda|\right)\le\frac{A}{| \Lambda|^{\alpha}}.
\end{equation*}
Note that condition (\ref{Halpha}) was exactly: there exists $\alpha>1$
such that $(H_\alpha)$ holds.
In order to have the coexistence result (Theorem \ref{coex}), we suppose  
moreover that $\Snu$ satisfies conditions (\ref{sup-des-moments}) and
$(H_\alpha)$ for some $\alpha>1$. 

\begin{theorem}
\label{temps-entier}
Let us suppose that $\Snu$ satisfies the previous general assumptions, and
assume moreover that:

1. The related semi-norm $\mu$ describing the directional asymptotic speeds is not identically null.

2. $\Snu$ is ``discrete'': there exists  a subset $\Den$
  of $\N$ such  
that $\Snu(\Den^{\Ed})=1$.

3. $\Snu$ satisfies the following finite energy property:  for each
  finite subset $\Lambda$ of $\Ed$ and 
each  $e_{\Lambda}\in \Den^{\Lambda}$, we have
$$\Snu(\omega_{\Lambda}=e_{\Lambda}|\mathcal{F}_{\Lambda^c})>0\quad
\Snu\as$$

4.
If $p=1$, we add the assumption: $\Den$ is unbounded.

5.
Some stronger integrability is assumed: suppose that one of the three
  following conditions is fulfilled
\begin{itemize}
\item
$(H_{\alpha})$ holds for some $\alpha>d^2+2d-1$.
\item
$p=1$ and the passage times of bonds have a moment of order $\alpha>d.$
\item
$p=1$, $\Snu$ is a product measure and the passage times of
   bonds have a second moment.
\end{itemize}
Then, for each pair $s_1,s_2$ of distinct sources in $\Z^d$,
$\P(\Coex(s_1,s_2))>0.$ \\
Moreover,  $$\P\left(\begin{array}{c}
    \mbox{there exists two disjoint semi-infinite geodesics}\\
    \mbox{starting from $0$ for the random distance $d$}
  \end{array}
\right)  >  0.$$
\end{theorem}

The last integrability condition is the only one that is specific to the 
discrete case: in the diffuse case of the previous section, we could give to a given edge an
arbitrary small value thanks to (\ref{energiefinie}), and there was no need
to control the length of an optimal path. Here, as passage times are
integers, we need a stronger integrability assumption that helps to control
these paths. The second moment assumption is classical in
\iid first-passage percolation to ensure the shape theorem --see the
reference article \cite{kesten} lemma 3.5; the $(H_\alpha)$ assumption with
$\alpha>d$ is the one used by Boivin in \cite{Boivin} for the shape
theorem in stationary first-passage percolation. Finally, the $(H_\alpha)$
assumption with $\alpha>d^2+2d-1$ is the one we use in
\cite{garet-marchand} lemma 3.7 to obtain the shape theorem when the edges
can be closed.
Note that, if $\Snu$ is the product measure $\nu^{\otimes\Ed}$,
assumption $(H_{\alpha})$ follows from the Marcinkiewicz-Zygmund inequality as soon as the passage time of an edge has a moment of order strictly
greater than $2\alpha$.

 In any case, we obtain the following estimate:

\begin{lemme}
\label{lem-controle}
There exists $K_1>0$ such that for every $a \in\Zd$, we can construct a random integer $M(a)<+\infty$ such that
$$a\communique\infty \text{ and }y\communique\infty\text{ and }\Vert y\Vert_1\ge M(a)\Longrightarrow d(a,y)\le K_1\Vert y\Vert_1.$$
\end{lemme}

\subsection*{Examples}

\begin{itemize}
\item Take  $p<1$ and $\Snu=\nu^{\otimes\Ed}$, where the support of $\nu$ is a
  finite subset of $\N^{*}$. As a special case, $\nu=\delta_1$ gives the
  classical chemical distance on a Bernoulli percolation cluster:
\begin{coro}[Geodesics on a Bernoulli cluster]
\label{geo-bernoulli}
For each $p>p_c$, consider Bernoulli percolation with parameter $p$.
Then, there almost surely exists a point of the infinite cluster
from which we can draw two disjoint semi-infinite geodesics.
\end{coro}

\begin{proof}
The considered event is translation-invariant, so its probability is null or full. By Theorem~\ref{temps-entier} with $\Snu=\delta_1^{\otimes\Ed}$, it can not be null. 
\end{proof}

\item Consider a Poisson point process on $\Rd$ with an intensity proportional to Lebesgue's measure. Let $\alpha\in\N^{*}$ and define the passage time $\eta_e$ by
$\eta_e=1+\alpha n_e$, where $n_e$ is the number of obstacles around $e$, \ie the number of points of the Poisson process which are closer from $e$ than from any other edge.

\end{itemize}

We can now begin the proof of Theorem \ref{temps-entier}.

\begin{proof} \textit{(Coexistence result).} The goal  is to prove that
  for each pair $s_1,s_2$ of distinct sources in $\Z^d$,
$\P(\Coex(s_1,s_2))>0.$ 

By translation invariance, we can suppose $s_1\in \Z^d
\backslash \{0\}$ and $s_2=0$. Since $\mu$ is not identically null, we can
find $x\in\Zd$ such that $\|s_1\|_1$ and $\|x\|_1$ have the same parity and
such that $\mu(x)\ne 0$. 
Thanks to Theorem \ref{coex}, we can consider an odd integer $n_0$ such
that $\P(\Coex(0,n_0 x))>0$. Note $s'_1=n_0 x$. We are going to prove that
$\P( \Coex(0,s_1))>0$.

Take $K_1>0$ and $M(0)$ and $M(s'_1)$ as defined in lemma \ref{lem-controle}. 
Since $$\lim_{n\to +\infty}\P(\{M(0) \le n \}\cap \{M(s'_1) \leq n \}\cap
\Coex(0,s'_1))=\P(\Coex(0,s'_1))>0,$$ we can find an integer $R_1$ such that
$$\P \left( \{M(0) \le R_1\}\cap \{M(s'_1) \leq R_1 \}\cap \Coex(0,s'_1)
\right)>0.$$ 

Let $\Lambda=\{x\in\Zd;\Vert x\Vert_1\le R\}$, for a large integer $R$ whose exact value will be fixed later.
 The idea is then to show that every
configuration  
$(\omega,\eta)$ in the event $A=\{M(0)\le R_1\}\cap \{M(s'_1)\le R_1\}\cap
\Coex(0,s'_1)$ can be modified inside 
the ball $\Lambda$ to get a configuration $(\omega',\eta')$ where
$\Coex(0,s_1)$  holds. A classical finite energy argument  concludes 
the proof: at first, note that $\P=\P_p\otimes\Snu$ also enjoys the finite energy property. Now if $B$ is a subset of $\Omega$ such that there exists
a map $f:A\to B$ with $f(x)_{\Lambda^c}=x_{\Lambda^c}$ for each $x\in A$,
then $P(B)>0$, because
\begin{eqnarray*}
\P(B) & = & \int_{\Omega}\P(B|\mathcal{F}_{\Lambda^c})(x)\ d\P(x)\\
      & \ge & \int_{A}\P(B|\mathcal{F}_{\Lambda^c})(x)\ d\P(x)\\
 & \ge & \int_{A}\P(\{f(x)\}|\mathcal{F}_{\Lambda^c})(x)\ d\P(x)\\
& > & 0.
\end{eqnarray*}

Let us explain now the modification inside $\Lambda$.
In the following, we will assume without loss of generality that the greatest
common divisor of the elements of $\Den$ is 1.
By the lemma of Bezout, we can find a finite family of integers $a_k$ and $s_k$, with $s_k\in \Den$, such that ${\sum_k} a_k s_k=1$.
Note 
$$
\begin{array}{ll}
\Den_{+}=\{k\in \Den; a_k>0\}, &  \Den_{-}=\{k\in \Den; a_k<0\}, \\
C_1=\sum_{k\in \Den_+} a_k, &  C_2=\sum_{k\in \Den_-} (-a_k), \\
b_1=\mbox{ smallest odd element of } \Den, & 
b_2=\mbox{ smallest even element of } \Den, \\
C=\max(C_1,C_2), & B=\max(b_1,b_2).  \\
\end{array}
$$
By convention, if $\Den$ only contains odd integers, we set $b_2=b_1$ and $B=0$.

The next lemma is a geometrical result, and we omit its proof because it is
rather tedious and not particularly illuminating:
\begin{lemme}
\label{un-beau-chemin}
Let us consider two fixed points $a_0,a_1\in\Zd$ (not necessarily distinct) and two non-negative numbers $D$ and $K$.
Let us  note $\Lambda_n=\{x\in\Zd;\Vert x\Vert_1\le n\}$.

There exists $\kappa=\kappa(a_0,a_1,D,K)<+\infty$ such that the following holds
as soon as $n\ge\kappa$:

For each distinct $r_0,r_1\in\Zd$ with $\Vert r_0\Vert_1=\Vert r_1\Vert_1=n$ and each
integer $l$ which has the same parity than $\Vert a_0-a_1\Vert_1$ and satisfies to
$|l|\le K n+D$ , one can
construct inside $\Lambda_n$  two simple paths $\gamma_0$ from $a_0$ to
$r_0$ and $\gamma_1$ from $a_1$ to $r_1$ with no common point (but maybe
$a_0$ if $a_0=a_1$) and such that:
$$|\gamma_0 | \geq K n+D, \; \; |\gamma_1 | \geq K n+D, \; \; |\gamma_0 | -|\gamma_1|=l.$$
 \end{lemme}
  
We can now define the radius $$R= \max(\kappa(0,s_1,B,K_1 C), \|s'_1\|_1+2,
R_1)$$ and define $\Lambda=\Lambda_R$.
Consider a semi-infinite geodesic starting from $0$ (\resp $s'_1$) and define by $r_0$ (\resp $r_1$) the 
last point of this semi-infinite geodesic which belongs to $\Lambda$.
Denote
$$L=d(r_0,0)-d(r_1,s'_1).$$ 
For simplicity, we will suppose, without loss of generality, that $L$ is
non-negative. Remember that $s'_1$ has the same parity than $s_1$.
Let us define 
$$
\left\{
\begin{array}{rl}
b'_2=b_2 \mbox{ and }b'_1=b_1 & \mbox{ if } \Vert s_1 \Vert_1 \mbox{ does
  not have the same parity than } L(C_1-C_2), \\
b'_1=b'_2=0 & \mbox{ otherwise}. 
\end{array}
\right.
$$
$$\text{ and }l=L(C_1-C_2)+b'_1-b'_2.$$

Note that $b_1-b_2$ is odd, unless $\Den$ only contains odd integers.
But in that case,  $C_1-C_2$ is odd and $L$ has the same parity than 
$\Vert r_0\Vert_1+\Vert r_1 \Vert_1+\Vert s'_1\Vert_1$, that is the same parity than $\Vert s_1\Vert_1$.

Thus, $l$ and $\Vert s_1\Vert_1$ always have the same parity. Note that 
\begin{eqnarray*}
|l| & \le & |L| |C_1-C_2|+B \\
& \le & \max(d(0,r_0),d(s'_1,r_1)).\max(C_1,C_2)+B \\
& \le & K_1 .\max(\Vert r_0\Vert_1,\Vert r_1\Vert_1).C+B=B+CK_1R. 
\end{eqnarray*}
So, by lemma~\ref{un-beau-chemin}, and by the choice we made for $R$, 
 one can
construct inside $\Lambda$ two simple paths with no common point
$\gamma_0$ from $0$ to $r_0$, and $\gamma_1$ from $s_1$ to $r_1$ such that 
$$|\gamma_0 | \geq CK_1R+B, \; \; |\gamma_1 | \geq CK_1R+B, \; \; |\gamma_0 | -|\gamma_1|=l.$$
Let us note $k=|\gamma_0|-(LC_1+b'_1)$. As proved in the upper bound for
$|l|$, $LC_1 \leq  CK_1R.$ We thus have 
$$k\ge CK_1R+B -(LC_1+b'_1)\ge B-b_1 \geq  0.$$
Obviously, $|\gamma_0|=LC_1+b'_1+k$ and
$|\gamma_1|=LC_2+b'_2+k$. Define also the following quantity $M$ that will
played the role of an ``infinite'' passage time for open edges:
\begin{eqnarray*}
M & = & \max \left\{ L\sum_{i\in \Den_+} a_i s_i + b'_1b_2+kb_1, 
 L\sum_{i\in \Den_-} (-a_i) s_i + (b'_2+k)b_1 \right\} \\
&& + \max \left\{ d(x,y)(0_\Lambda\omega_{\Lambda^c}, \eta),
  \|x\|_1=\|y\|_1=R \right\}.
\end{eqnarray*}
Note that $M$ is in $\mathcal F_{\Lambda^c}$, the $\sigma$-algebra
generated by  $\{(\omega_e,\eta_e), e \in \Lambda^c\}$. 
Now define, for every $(\omega_e,\eta_e) \in A$,  the configuration
$(\omega',\eta')\in \Omega$: set 
$(\omega'_e,\eta'_e)=(\omega_e,\eta_e)$ for $e\in\Ed\backslash\Lambda$ 
and define
$(\omega'_{\Lambda},\eta'_{\Lambda})$ inside $\Lambda$ as follows:

\begin{enumerate}
\item If $p<1$, $\forall e \in \Lambda\backslash (\gamma_0 \cup \gamma_1),
  \omega'_e=0$ and $\eta'_e=b_1$, but this value does not play a special role; \\
if $p=1$, $\forall e \in \Lambda\backslash (\gamma_0 \cup \gamma_1),
\eta'_e>M$, and $\omega'_e=1$. 
\item $\forall e \in \gamma_0 \cup \gamma_1, \omega'_e=1$.
\item Assign a passage time to edges in $\gamma_0$ as follows (remember
  that $|\gamma_0|=LC_1+b'_1+k$):
first, for each $i\in \Den_+$, give to $a_iL$ edges the value $\eta'_e=s_i$
and next complete giving to $b'_1$ other edges the value $\eta'_e=b_2$ and
to $k$ other edges the value $\eta'_e=b_1$.
\item Assign a passage time to edges in $\gamma_1$ as follows (remember
  that $|\gamma_1|=LC_2+b'_2+k$):
first, for each $i\in \Den_-$, give to $-a_iL$ edges the value $\eta'_e=s_i$
and next complete giving to the remaining $b'_2+k$ edges the value $\eta'_e=b_1$.
\end{enumerate}
Now we immediately obtain: 
\begin{eqnarray*}
 \sum_{e\in\gamma_0} \eta'_e & = & L\sum_{i\in \Den_+} a_i s_i+b'_1b_2+kb_1,
\\ 
 \sum_{e\in\gamma_1} \eta'_e & = & L\sum_{i\in \Den_-} (-a_i)
s_i+(b'_2+k)b_1, \\
 \sum_{e\in\gamma_0} \eta'_e-\sum_{e\in\gamma_1} \eta'_e & = & L\big(
\sum_{i\in \Den} a_i s_i\big)+ b'_1b_2-b'_2b_1=L= d(0,r_0)-d(s'_1,r_1).
\end{eqnarray*}

For $(\omega,\eta) \in A \subset \Coex(0,s'_1)$, there exist two infinite
geodesics starting from $0$ and $s'_1$. Let us denote by $\Gamma_0$ (\resp
$\Gamma_1$) the part beginning at $r_0$ 
(\resp $r_1$) in the geodesic  starting
from $0$  (\resp $s'_1$) in the configuration $(\omega,\eta)$.

We are going to prove that $\gamma_0\cup\Gamma_0$ (\resp $\gamma_1\cup\Gamma_1$) is an infinite geodesic starting from $0$ (\resp $s_1$) in the configuration $(\omega',\eta')$.
Let $x$ be a point of $\Gamma_0$. Let us prove that an optimal path from $0$
to $x$ in the configuration $(\omega',\eta')$ is included in $\gamma_0 \cup
\Gamma_0$.  

Let $\gamma$ be an optimal path from $0$ to $x$ in 
the configuration $(\omega',\eta')$, and denote by $z$ the point from which
the path $\gamma$ exits from $\Lambda$. We have 
$$
d(0,z)=\sum_{e \in \gamma} \eta'(e) \leq \sum_{e \in \gamma_0}
\eta'(e)+d(r_0,z).
$$
But since $\displaystyle\sum_{e \in \gamma_0} \eta'(e)=L\sum_{i\in \Den_+} a_i
s_i+b'_1b_2+kb_1$ and $d(r_0,z)(\omega',\eta') \leq
d(r_0,z)(0_\Lambda\omega_{\Lambda^c},\eta),$ it follows that $d(0,z) \leq
M$. By definition of $M$, it ensures that $\gamma$ does not use any bond in
$\Lambda$, except those used in $\gamma_0 \cup \gamma_1$, and particularly,
it implies that $z=r_0$, and thus an optimal path from $0$ to $x$ is
included in $\gamma_0 \cup \Gamma_0$. 

Similarly, let $\gamma$ be an optimal path from $s_1$ to $x$ in 
the configuration $(\omega',\eta')$, and denote by $z$ the point from which
the path $\gamma$ exits from $\Lambda$. We have 
$$
d(0,z)=\sum_{e \in \gamma} \eta'(e) \leq \sum_{e \in \gamma_1}
\eta'(e)+d(r_1,z).
$$
But since $\sum_{e \in \gamma_1} \eta'(e)=L\sum_{i\in \Den_-} (-a_i)
s_i+(b'_2+k)b_1$ and $d(r_1,z)(\omega',\eta') \leq
d(r_1,z)(0_\Lambda\omega_{\Lambda^c},\eta),$ it follows that $d(s_1,z) \leq
M$. By definition of $M$, it ensures that $\gamma$ do not use any bond in
$\Lambda$, except those used in $\gamma_0 \cup \gamma_1$, and particularly,
it implies that $z=r_1$, and thus an optimal path from $s_1$ to $x$ uses
$\gamma_1$ to exit from $\Lambda$. 

Let us now prove that if $x \in \Gamma_0$, $d(0,x)(\omega',\eta')<d(s_1,x)(\omega',\eta')$. 
\begin{eqnarray*}
d(s_1,x)(\omega',\eta') 
& = & d(s_1,r_1)(\omega',\eta')+d(r_1,x)(0_\Lambda \omega_{\Lambda^c},\eta'), \\
& = & d(s_1,r_1)(\omega',\eta')+d(r_1,x)(0_\Lambda \omega_{\Lambda^c},\eta), \\
d(0,x)(\omega',\eta')
& = & d(0,r_0)(\omega',\eta') +d(r_0,x)(0_\Lambda \omega_{\Lambda^c},\eta'), \\
& = & d(0,r_0)(\omega',\eta') +d(r_0,x)(0_\Lambda \omega_{\Lambda^c},\eta).
\end{eqnarray*}
Consequently, 
\begin{eqnarray*}
&& (d(s_1,x)-d(0,x))(\omega',\eta') \\
&  = & (d(s_1,r_1)-d(0,r_0))(\omega',\eta')
+ d(r_1,x)(0_\Lambda \omega_{\Lambda^c},\eta)
- d(r_0,x)(0_\Lambda \omega_{\Lambda^c},\eta), \\
& = & (d(s'_1,r_1)-d(0,r_0)) (\omega,\eta) 
+ d(r_1,x)(0_\Lambda \omega_{\Lambda^c},\eta)
- d(r_0,x)(0_\Lambda \omega_{\Lambda^c},\eta), \\
& \geq  & d(s'_1,x) (\omega,\eta) - d(0, x)(\omega,\eta)>0,
\end{eqnarray*}
because $x \in \Gamma_0$, which is a part of the infinite geodesic issued
from $0$ in the configuration $(\omega, \eta)$. Thus $\gamma_0 \cup
\Gamma_0$ is an infinite geodesic issued from $0$ in the configuration
$(\omega', \eta')$. 

In the same manner, working symmetrically with $\Gamma_1$, we prove that
$\gamma_1 \cup \Gamma_1$ is an infinite geodesic issued from $s_1$ in the
configuration 
$(\omega', \eta')$, and finally $(\omega', \eta') \in \Coex(0,s_1)$.
\end{proof}

\begin{proof} \textit{(Geodesics result).} The goal here is to prove that
$$\P\left(\begin{array}{c}
    \mbox{there exists two distinct semi-infinite geodesics}\\
    \mbox{starting from $0$ for the random distance $d$}
  \end{array}
\right)  >  0.$$
The proof is exactly the same as the previous one. The only difference is
to use the single source $0$ rather than two distinct sources $0, s'_1$. The geometrical
structure of the modification is once again given by lemma
\ref{un-beau-chemin}, and the adjustment of the values is made as before.
\end{proof}

As seen previously, a trouble with integer passage times is that some
points can be reached at the very same moment by  the two distinct
infections. This case can be ruled out under some extra assumptions, and
this is the goal of the next result. But first, 
for two distinct sources $x,y \in \Z^d$, we say that the event
$\SepCoex(x,y)$ happens if 
$$\left\{
\begin{array}{ll}
& \{z\in \Z^d, d(x,z)<d(y,z)\}\text{ is infinite}\\
\text{ and } & \{z\in \Z^d, d(x,z)>d(y,z)\}\text{ is infinite} \\
\text{ and } & \forall z \in \Z^d, \; d(x,z) \neq d(y,z)\text{ unless }d(x,z)=d(y,z)=+\infty.
\end{array}
\right.$$
We have the following result: 

\begin{lemme}
\label{coro-impair}
Denote by $\mathcal O$ the set of non-negative odd integers, and as previously, 
let $d \geq 2$, $p>p_c(d)$, $\Snu$ a
  stationary ergodic probability measure on $\mathcal O^{\Ed}$
  satisfying (\ref{sup-des-moments}) and (\ref{Halpha}). 
Then,  for $x\in\Zd$ with $\Vert x\Vert_1$ odd,
$$\P(Coex(0,x)\backslash\SepCoex(0,x))=0.$$
\end{lemme}

\begin{proof}

By the assumption we made on $\mu$, the length of a path
from $x$ to $y$ has the same parity than $\Vert x-y\Vert_1$.
So, the identity $d(x,z)=d(y,z)$ can only happens if $\Vert x-y\Vert_1$
is even.
\end{proof}

Now, for a given point $x$ with $\Vert x\Vert_1$ odd, 
the fact that $\P(\SepCoex(0,x))>0$ can  be obtained as a consequence
of Theorem~\ref{coex} or Theorem~\ref{temps-entier}.
Note also that when the assumptions of lemma~\ref{coro-impair} are fulfilled,
$\mu$ is always a norm: since the passage time of a bond is an odd integer, it is at least equal to 1. Then, it is easy to see that for each $x\in\Zd$, we have $\mu(x)\ge \|x\|_1$.


\section{An example of a discrete time competing process}


The last section is devoted to the study of a natural example
of a non-trivial dynamical system which can be studied with the help of
Theorem~\ref{temps-entier} and lemma~\ref{coro-impair}.

Consider two species, say blue and yellow, which attempt to conquer the space $\Zd$. At each instant, each fertile cell tries to contaminate each of its non-occupied neighbors. It succeeds with probability $p$. In case of success, the non-occupied cell takes the color of the infector. If a yellow cell and blue cell simultaneously succeed in contaminating a given cell, this one takes the green color. If a green cell and another cell simultaneously succeed in contaminating a given cell, this one takes the green color. At the next step, the individuals that have just been generated are fertile, but the previous generation is no more fertile. 
We make the following assumptions:
\begin{itemize}
\item  the success of each attempt of contamination at a given time does not depend on the past,
\item the successes of simultaneous attempts to contamination are independent.
\end{itemize}

The first assumption allow a modelization by an homogeneous  Markov chain.
Markov chains satisfying to the second condition are sometimes called Probabilistic Cellular Automata (PCA).
 
Let us define 
$$S=\{0,\blue,\yellow,\green,\bluestar,\yellowstar,\greenstar\},$$
where $0$ is the state of an empty cell, $\blue,\yellow,\green$ the states of  fertile cells, and $\bluestar,\yellowstar,\greenstar$ the states of unfertile cells.

Since we will study the evolution of a system which starts with only two cells, we will only deal with configurations in which a finite numbers of cells
are nonempty. So, we will deal with a classical Markov chain on the denumerable set
$$C=\{\xi\in S^{\Zd};\quad\exists \Lambda\text{ finite},\quad \xi_k=0\text{ for }k\in\Zd\backslash\Lambda\}.$$

We now define for $\text{color}\in \Act=\{\blue,\yellow,\green\}:$
$$n(\text{color},x)(\xi)=|\{y\in\Zd; \|x-y\|_1=1\text{ and }\xi_y=\text{color}\}|,$$
and
$$s(\text{color},x)=1-(1-p)^{n(\text{color},x)},$$
which represents the probability that at least one neighbor of $x$ succeeds in  infecting $x$ with the given  color. 

The considered dynamics form an homogeneous  PCA with space state $S$
 and whose local evolution rules are given by  
$$
p_x(s,t)  = \left\{ 
\begin{array}{l}
\begin{array}{ll}
(1-p)^{n(\yellow,x)+n(\blue,x)+n(\green,x)} & \text{if }s=t=0\\
 s(\blue,x) (1-p)^{n(\yellow,x)+n(\green,x)} & \text{if }s=0\text{ and }t=\blue\\
 s(\yellow,x)(1-p)^{n(\blue,x)+n(\green,x)} & \text{if }s=0\text{ and }t=\yellow\\
s(\green,x)+(1-s(\green,x))s(\blue,x)s(\yellow,x)& \text{if }s=0\text{ and }t=\green\\ 
\end{array}\\
\begin{array}{ll}
1 & \quad\quad\quad\quad\quad\quad\quad\quad\quad\quad\quad\quad\text{if }s\in\{\blue,\yellow,\green\}\text{ and }t=s^*\\
1 & \quad\quad\quad\quad\quad\quad\quad\quad\quad\quad\quad\quad\text{if }s\in\{\blue^*,\yellow^*,\green^*\}\text{ and }t=s\\
0 & \quad\quad\quad\quad\quad\quad\quad\quad\quad\quad\quad\quad\text{otherwise}.
\end{array}
\end{array}
\right.$$
In term of Markov chains, it means that the transition matrix is defined by
$$\forall (\xi,\omega)\in C\times C\quad p(\xi,\omega)=\prod_{k\in\Zd} p_k(\xi_k,\omega_k).$$
The product is convergent because only a finite numbers of terms differs from $1$.

With the help of the tools that we have developed above, we will prove
the following theorem:

\begin{theorem}
Let $p>p_c$. For $\sy,\sb\in\Zd$ with $\sb\ne\sy$, let us denote by
$\P_{p,\sy,\sb}$ the law of a PCA $(X_n)_{n\ge 0}$ following the dynamics
described above, and  starting a configuration with exactly two non-empty cells:
a blue cell at site $\sb$, a yellow cell at site $\sy$, the others cells being empty. Then, 
$$ \P_{p,\sy,\sb}(\forall n\in\N\;\exists (x,y)\in\Zd\times\Zd,\; X_n(x)=\blue\text{ and } X_n(y)=\yellow)>0.$$
If moreover   
$\|\sy-\sb\|_1$
is odd, 
green cells never appear.
\end{theorem}

\begin{figure}
\begin{tabular}{lcr}
\includegraphics{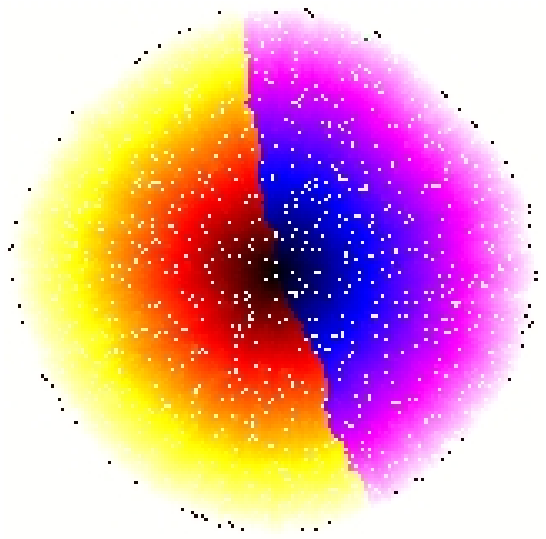} & & \includegraphics{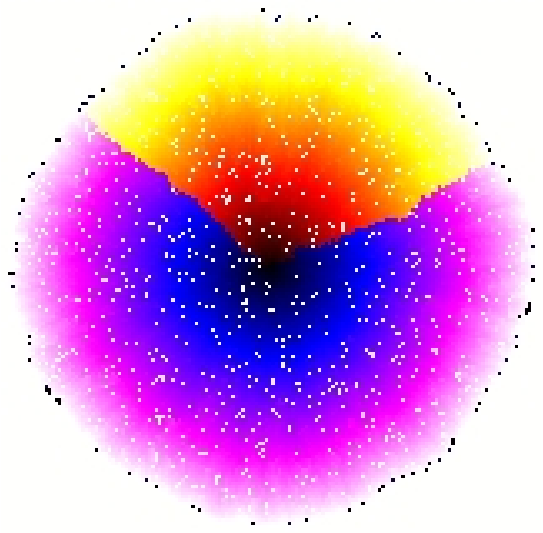}
\end{tabular}
\caption{Two samples of simulation of the competing process when $p=0.6$. The
process is stopped when the border of a given box is attained by one of the two species. The color in the picture is determined by the time of coloring and the type of the cell.}
\end{figure}

The following lemma gives the link between this PCA and our competing model.

\begin{lemme}
\label{lien}
Consider a probability space where lives a family $(\omega_e)_{e\in\Ed}$ of independent Bernoulli variables with parameter $p$, which defines a  random chemical distance $D$.

Let $\sy,\sb\in\Zd$ with $\sb\ne\sy$. Define
\begin{eqnarray*}
X_n(x) & = &
\begin{cases}
 \blue & \text{if }n=D(\sb,x)<D(\sy,x)\\
 \yellow & \text{if }n=D(\sy,x)<D(\sb,x)\\
 \green & \text{if }n=D(\sy,x)=D(\sb,x)\\
 \bluestar & \text{if }D(\sb,x)<\min(D(\sy,x),n)\\
 \yellowstar & \text{if }D(\sy,x)<\min(D(\sb,x),n)\\
 \greenstar & \text{if }D(\sy,x)=D(\sb,x)<n\\
 0 & \text{otherwise}.
\end{cases}
\end{eqnarray*}
Then, $(X_n)_{n\ge 0}$ is an homogeneous  PCA with space state
$$S=\{0,\blue,\yellow,\green,\bluestar,\yellowstar,\greenstar\}$$ associated to the probabilities $p_x(s,t)$ defined above.
\end{lemme}

\begin{proof}
Let us consider the map 
\begin{eqnarray*}
f:S^{\Zd}\times\Omega_E & \to & S^{\Zd}\\
(\xi,\omega) & \mapsto &(f_x(\xi_x,\omega))_{x\in\Zd}
\end{eqnarray*}
where $f_x:S\times\Omega_E\to S$ is defined by
$$\left\{
\begin{array}{lll}
f_x(s,\omega) & = & s^{*} \text{ for each }s\in \{\blue,\yellow,\green\}\\
f_x(s,\omega) & = & s \text{ for each }s\in \{\blue^*,\yellow^*,\green^*\}\\
f_x(0,t) & = & 
\begin{cases}
\blue & \text{if } \Act\cap\{\xi_y; \|x-y\|_1=1\text{ and }\omega_{\{x,y\}}=1\}=\{\blue\}\\
\yellow & \text{if } \Act\cap\{\xi_y; \|x-y\|_1=1\text{ and }\omega_{\{x,y\}}=1\}=\{\yellow\}\\
\green & \text{if } \Act\cap\{\xi_y; \|x-y\|_1=1\text{ and }\omega_{\{x,y\}}=1\}=\{\blue,\yellow\}\\
\green & \text{if } \Act\cap\{\xi_y; \|x-y\|_1=1\text{ and }\omega_{\{x,y\}}=1\}\supset\{\green\}\\
0 & \text{otherwise}.
\end{cases}
\end{array}
\right.$$
By considering Dijkstra's algorithm in the particular case where the travel times are constant, it is not difficult to see that $(X_n)_{n\ge 0}$ satisfy to 
the recurrence formula $X_{n+1}=f(X_n,\omega)$.
To recognize $(X_{n})_{n\ge 0}$ as a convenient PCA, we will build a coupling
of $\omega$ with a i.i.d. sequence $(\omega^n)_{n\ge 1}$ to obtain the canonical Markov Chain representation   $X_{n+1}=f(X_n,\omega^n)$.

Let $(\Omega,\mathcal{F},P)$ be a probability space with $\zeta^0,\omega^0,\omega^1,\omega^2,\dots$ independent $\{0,1\}^{\Ed}$ valued variables with $\Ber(p)^{\otimes\Ed}$ as common law.

We define $A_0=\{\sb,\sy\}$ and recursively
\begin{eqnarray*}&
\begin{cases}
B_{n+1}  = & \{y\in \Zd\backslash A_n\quad \exists x\in\partial A_n: \|x-y\|_1=1\text{ and }\omega_{\{x,y\}}^n=1\}\\
A_{n+1} = & A_n\cup B_{n+1}.
\end{cases}
\end{eqnarray*}
Note that the random set $B_{n+1}$ is measurable with respect to the $\sigma$-algebra generated by $(\omega^0,\omega^1,\dots,\omega^n)$.
We define $\zeta^n$  recursively by
\begin{eqnarray*}
\zeta_e^{n+1} = &
\begin{cases}
\omega_e^{n+1} & \text{if }e=\{x,y\}\text{ with } (x,y)\in \partial A_n\times \Zd\backslash A_{n}\\
\zeta_e^n  & \text{otherwise.}
\end{cases}
\end{eqnarray*}
By natural induction, we prove that the law of $\zeta^n$ under $P$ is $\Ber(p)^{\otimes\Ed}$.
By construction, each bond $e$ writes  $e=\{x,y\}\text{ with } (x,y)\in \partial A_n\times \Zd\backslash A_{n}$ for at most one value of $n$.
It follows that the sequence $\zeta^n$ converges in the product topology. 
Let us denote by $\omega^{\infty}$ its limit. Since the  law of $\zeta^n$ under $P$ is $\Ber(p)^{\otimes\Ed}$, it follows that the law of $\omega^{\infty}$
 under $P$ is also $\Ber(p)^{\otimes\Ed}$.

Now, it is not difficult to see that sequence $(X_n)_{n\ge 0}$ defined from 
 $\omega^{\infty}$ as previously, satisfies to the recurrence 
formula $X_{n+1}=f(X_n,\omega^{\infty})$, but also to $X_{n+1}=f(X_n,\omega^n)$.

It is now proved that $(X_{n})_{n\ge 0}$ is an homogeneous Markov chain.
The recognition of the transition matrix follows from an elementary calculus.
\end{proof}

We can now prove the theorem announced above.

\begin{proof}
Clearly, lemma~\ref{lien} connects the considered PCA with the random distance
studied in Theorem~\ref{temps-entier}.
Here, the passage times of open bonds are identically equal to $1$, which
is obviously an odd number. 
By lemma~\ref{coro-impair}, this prevents from the appearance of green cells when $\|\sy-\sb\|_1$ is odd.
\end{proof}

\subsection*{Remarks.}
If $\|\sy-\sb\|_1\ne 0$ is even and if the two species infinitely grow, there
are necessarily green cells at the boundary between blue cells and yellow cells.
   
\begin{figure}
\label{chantdudepart}
\includegraphics[scale=0.3]{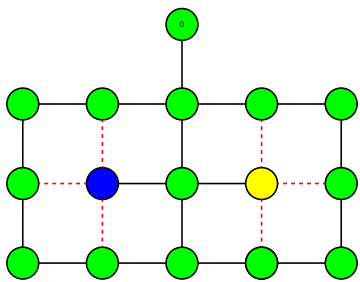} 
\caption{Green surrounding the two sources}
\end{figure}

A natural question is the following: is it possible to have an
infinite set of green cells surrounding the blue cells and the yellow cells ? The answer is yes, as soon as $\|\sy-\sb\|_1\ne 0$ is even:
consider figure~\ref{chantdudepart}.

The picture describes a particular case when $d=2$, but the reasoning
can obviously be generalized.

In this case, the yellow flow and the blue flow immediately converge
to engender a green flow. They also do not develop themselves elsewhere.
If the point labelled~$0$ belongs to the infinite cluster, then the result is 
proved.

It is now to see that, conditionally to the states of the bonds imposed by this
picture, the probability that  $0$ belongs to the infinite cluster is strictly positive, which follows from a classical modification argument.


\bibliographystyle{plain}
\bibliography{competition-p}

\begin{thebibliography}{1}

\bibitem{Boivin}
Daniel Boivin.
\newblock First passage percolation: the stationary case.
\newblock {\em Probab. Theory Related Fields}, 86(4):491--499, 1990.

\bibitem{Deijfen-Haggstrom}
Maria Deijfen and Olle Häggström.
\newblock The initial configuration is irrelevant for the possibility of mutual
  unbounded growth in the two-type richardson model.
\newblock {\em Preprint, available at \texttt{http://www.math.su.se/\~{
  }mia/sorm.pdf}}, 2003.

\bibitem{garet-marchand}
Olivier Garet and Régine Marchand.
\newblock Asymptotic shape for the chemical distance and first-passage
  percolation in random environment.
\newblock {\em Preprint, available at
  \texttt{http://arxiv.org/abs/math.PR/0304144}}, 2003.

\bibitem{grimmett}
Geoffrey Grimmett.
\newblock {\em Percolation}, volume 321 of {\em Grundlehren der Mathematischen
  Wissenschaften [Fundamental Principles of Mathematical Sciences]}.
\newblock Springer-Verlag, Berlin, second edition, 1999.

\bibitem{Haggstrom-Pemantle-1}
Olle H{\"a}ggstr{\"o}m and Robin Pemantle.
\newblock First passage percolation and a model for competing spatial growth.
\newblock {\em J. Appl. Probab.}, 35(3):683--692, 1998.

\bibitem{Harris}
T.~E. Harris.
\newblock A lower bound for the critical probability in a certain percolation
  process.
\newblock {\em Proc. Cambridge Philos. Soc.}, 56:13--20, 1960.

\bibitem{kesten}
Harry Kesten.
\newblock Aspects of first passage percolation.
\newblock In {\em \'Ecole d'\'et\'e de probabilit\'es de Saint-Flour,
  XIV---1984}, volume 1180 of {\em Lecture Notes in Math.}, pages 125--264.
  Springer, Berlin, 1986.

\bibitem{stout}
William~F. Stout.
\newblock {\em Almost sure convergence}.
\newblock Academic Press [A subsidiary of Harcourt Brace Jovanovich,
  Publishers], New York-London, 1974.
\newblock Probability and Mathematical Statistics, Vol. 24.

\end{thebibliography}

\end{document}